\documentclass[leqno,12pt]{article}
\usepackage{amsmath, amsthm}
\usepackage{amsfonts}
\usepackage{amssymb}
\usepackage{graphicx}
\usepackage{color}
\usepackage{hyperref}

\newtheorem{theorem}{Theorem}[section]
\newtheorem{corollary}[theorem]{Corollary}
\newtheorem{lemma}[theorem]{Lemma}
\newtheorem{proposition}{Proposition}[section]

\definecolor{blue}{rgb}{0.05, .05, .65}

\newcommand{\re}{\mathbb{R}}
\newcommand{\ren}{\mathbb{R}^N}
\newcommand{\ve}{\varepsilon}

\newcommand{\dx}{\,{\rm d}x}
\newcommand{\dy}{\,{\rm d}y}

\numberwithin{equation}{section}
\def\qed{\,\unskip\kern 6pt \penalty 500
\raise -2pt\hbox{\vrule \vbox to8pt{\hrule width 6pt
\vfill\hrule}\vrule}\par}

\newcommand{\nc}{\normalcolor}

\title{\bf The evolution fractional p-Laplacian equation in $\ren$. Fundamental solution \\ and asymptotic behaviour}
\author{ \sc \large J. L. V\'azquez\footnote{Univ. Aut\'onoma de Madrid, Spain.}}

\date{ }
\begin{document}
\maketitle

\begin{abstract} We consider the natural time-dependent fractional $p$-Laplacian equation posed in the whole Euclidean space, with parameters $p>2$ and $s\in (0,1)$ (fractional exponent).  We show that the  Cauchy Problem for data in the Lebesgue $L^q$ spaces is well posed, and show that the solutions form a family of  non-expansive semigroups with regularity and other interesting properties. As  main results, we construct the self-similar fundamental solution for every mass value $M,$ and prove that general finite-mass solutions converge  towards that fundamental solution having the same mass, and convergence holds in all $L^q$ spaces. A number of additional properties and estimates complete the picture.

\end{abstract}

\section{Introduction. The problem}

The nonlocal energy functional
\begin{equation}\label{Jsp1}
{\mathcal J}_{p,s}(u)= \frac1{p}\int_{\ren }\int_{\ren } \frac{|u(x)-u(y)|^p}{|x-y|^{N+sp}}\,dxdy\,.
\end{equation}
is a power-like functional with nonlocal kernel of the $s$-Laplacian type that has attracted a great deal of attention in recent years. It is just the $p$-power of the Gagliardo seminorm, used in the definition of the $W^{s,p}$ spaces (fractional Sobolev, Slobodeckii or Gagliardo spaces) with seminorm and norm given by
 $$
 \quad [u]_{s,p}^p=p{\mathcal J}_{p,s}(u), \qquad \|u\|_{s,p}^p=\int |u|^p\,dx + p{\mathcal J}_{p,s}(u),
 $$
 cf. \cite{AdFou, DiNPV}. We  consider the functional   for   exponents $0<s<1$ and $1<p<\infty$ in dimensions $N\ge 1$. Its  subdifferential 
${\mathcal L}_{s,p}$  is the nonlinear operator defined  a.e. by the formula
\begin{equation}\label{frplap.op}
\displaystyle \qquad {\mathcal L}_{s,p}(u):= P.V.\int_{\ren}\frac{\Phi(u(x,t)-u(y,t))}{|x-y|^{N+sp}}\,dy\,,
\end{equation}
where  we write $\Phi(z)=|z|^{p-2}z.$ It is a usually called the $s$-fractional $p$-Laplacian operator. It is then well-known from general theory that ${\mathcal L}_{s,p}$ is a maximal monotone operator in $L^2(\ren)$ with dense domain.

In this paper we study the corresponding gradient flow, i.e., the evolution equation
\begin{equation}\label{frplap.eq}
\partial_t u + {\mathcal L}_{s,p} u=0
\end{equation}
posed in the Euclidean space $x\in \ren$, $ N\ge 1$, for $t>0$. We will often refer to it as the EFPL equation (evolution fractional $p$-Laplacian equation). We supplement the equation with an initial datum
\begin{equation}\label{init.data}
\lim_{t\to 0} u(x,t)= u_0(x),
\end{equation}
where in principle $u_0\in L^2(\ren)$. However, the theory shows that Equation \eqref{frplap.eq} generates a continuous nonlinear semigroup in any $L^q(\ren)$ space, $1\le q<\infty$, in fact it is a nonexpansive semigroup for every $s,p$ and $q$ as specified.

We will concentrate in this paper on data $u_0\in L^1(\ren)$, which leads to the class of finite-mass solutions and  produces a specially rich theory. We will take all fractional exponents $0<s<1$, but restrict $p$ to the superlinear range $p>2$ for convenience of the study to be developed  here. The case $1<p<2$ is worth its own study. We will be specially interested in taking a Dirac delta as initial datum. In that case the solution is called a {\sl fundamental solution, \rm and also a {\sl source-type solution} (mainly in the Russian literature). Actually, there is one-parameter family of such solutions with parameter the mass $M$ of the initial point distribution. These fundamental solutions turn out to be asymptotic attractors for all integrable solutions in the precise sense to be stated in a moment.

The nonlinear fractional operator we are dealing with in this paper was mentioned in the paper \cite{IshiiN} by Ishii and Nakamura, see also \cite{ChLM12} by Chambolle, Lindgren and Monneau. There are  a number of works that cover the  evolution equation \eqref{frplap.eq} in the case where the space domain is a bounded subdomain $\Omega\subset \ren$, see \cite{MazRT, Vaz16, GalWar17} and references. References to the equation posed in the whole space are more recent, like \cite{BLS2019}, \cite{Strom2019}. Let us point out that inserting a constant factor in the definition of operator
 ${\mathcal L}_{s,p}$ does not change the properties of the solutions of Problem \eqref{frplap.eq}-\eqref{init.data} but it will be important in Section \ref{sec.limitcases}.

\subsection{Motivation and related equations}
Before presenting our results, let us briefly comment on motivations and closely  related equations.
Nonlinear equations like \eqref{frplap.eq} are motivated by the current interest in studying the combination of nonlinear and nonlocal terms in the formulation of the basic models of nonlinear diffusion in view of a large number of applications. These arise  in fields like  anomalous transport and diffusion, stochastic processes,  finance, elasticity,  conservation laws, porous medium flow, quasigeostrophic flows, image processing, population dynamics, flame propagation, chemical reactions of liquids, and so on.  It also has a theoretical interest for PDEs, Nonlinear Functional Analysis and Potential Theory. Some of these nonlinear nonlocal diffusion models are presented in the survey paper \cite{Vaz16}, where the nonlinearities are mainly of porous medium type, see also \cite{VazAbel, Vaz2014}.

The simplest equation in the fractional family is found in the limit case where $p=2$
\begin{equation}\label{frac.heq}
u_t+(-\Delta)^{s} u=0\,,
\end{equation}
i.\,e., the  heat equation associated to the  fractional Laplacian $(-\Delta)^{s}$, a nonlocal generalization of the Laplace operator that has been studied in classical monographs like \cite{Landkof, Stein}. The $s$-Laplacian is a linear operator that coincides with ${\mathcal L}_{s,2}$ up to a constant. Equation \eqref{frac.heq} inherits many of the well-known properties of the classical heat equation (case $s=1$) except for rates of space propagation, reflected in the fact the solutions with compactly supported data develop, for all positive times, spatial profiles with tails at infinity that decay like a power of distance, $u(x,t)\sim c(t)|x|^{-(N+2s)}$.  The equation has been amply discussed in the literature, see recent results in \cite{BPSV, BSV17, GarFT} on the existence and regularity theory, and \cite{Va18} for the asymptotic behaviour.

On the other hand, it is proved that in the limit $s\to1$ with $p\ne 2$, and after inserting a normalizing constant, we get the well-known evolution $p$-Laplacian  equation \
$$
\partial_t u=\Delta_p(u):=\nabla\cdot(|\nabla u|^{p-2}\nabla u),
$$
that has also been widely investigated since the early 1970's because of a number of applications (cf. for instance \cite{VazSmooth}, Section 11) and for its remarkable mathematical properties.
The semigroup method proved to be an effective method to treat the equation, see early works by B{\'e}nilan and V{\'e}ron \cite{Ben-pl,Veronsem}. Regularity theory is contained in the monograph by DiBenedetto \cite{DiBeBk}.
The recent literature on this topic is very large and has many novel features.


\subsection{Outline of the paper and main results}

 We focus on Problem \eqref{frplap.eq}-\eqref{init.data}, posed in $\ren.$ It is not difficult to prove that this Cauchy problem is well-posed in all $L^q(\ren)$ spaces, $1\le q< \infty$.
This parallels what is known in the case of bounded domains. In view of such works, we first review the main facts of the theory in Section \ref{sec.basic}. In particular, we  define the class of continuous strong solutions that correspond to $L^2$ and $L^1$  initial data and derive its main properties in detail.

We want to stress the differences brought about but the consideration of the whole space. In this respect, a most interesting question is that of finding the fundamental solution, i.e., the solution such that
\begin{equation}
\lim_{t\to 0} \int_{\ren} u(x,t)\varphi(x)\,dx=M\varphi(0),
\end{equation}
for every smooth and  compactly supported test function $\varphi\ge 0$, and some $M>0$.

\begin{theorem}\label{thm.ssfs} For every given mass $M>0$ there exists a unique self-similar solution of Problem \eqref{frplap.eq}-\eqref{init.data} with initial data $M\delta(x)$. It has the form
\begin{equation}\label{eq.ssf1}
U(x,t;M)=M^{sp\beta}t^{-\alpha}F(M^{-(p-2)\beta} x\,t^{-\beta})\,,
\end{equation}
with self-similarity exponents
\begin{equation}\label{eq.sse1}
\alpha=\beta N, \quad \beta=\frac{1}{N(p-2)+sp}\,.
\end{equation}
The profile $F(r)$  is a continuous, positive, radially symmetric ($r=|x|\,t^{-\beta}$), and decreasing function  such that $F(r) \approx r^{-(N+sp)}$ as $r\to\infty$.
\end{theorem}
 As usual, the sign $\approx$ means equivalence up to a constant factor, i.\,e., $0<c_1\le  F(r)r^{N+sp}\le c_2$. Here, $c_1$ and $c_2$ depend only on $N,s$ and $p$. We see that all fundamental solutions with $M>0$ are obtained from the one with unit mass, $M=1$, by a simple rescaling. Fundamental solutions with $M<0$ are obtained by just reversing the sign of the solution. For $M=0$ the fundamental solution becomes the null function. The theorem is proved in Section \ref{sec.fs}.  Important preliminaries take up Section
\ref{sec.barr}, where a delicate upper barrier is constructed; Section \ref{sec.masscon} about mass conservation; and Section \ref{sec.diff} about dissipation. The tail behaviour is settled in Section \ref{sec.pos.tail}, see in particular Corollary \ref{cor.ob.fs}. A numerical computation of the profile $F$ for different values of $s$ and $p$ is exhibited in Figures 1 and 2.

The fundamental solution is the key to the study of  the long-time behaviour of our problem with general initial data, since it represents,  in Barenblatt's words, the \sl intermediate asymptotics, cf. \rm \cite{Barbk96}. This is the  asymptotic result we obtain.

\begin{theorem}\label{thm.ab1} Let $u$ be a solution of Problem \eqref{frplap.eq}-\eqref{init.data} with initial data $u_0\in L^1(\ren)$ of integral $M$,  and let $U_M$ be the fundamental solution with that mass. Then,
\begin{equation}\label{lim.ab.L1}
\lim_{t\to\infty} \|u(t)-U_M(t)\|_1=0\,.
\end{equation}
We also have the $L^\infty$-estimate
\begin{equation}\label{lim.ab.inf}
\lim_{t\to\infty} t^{\alpha}\|u(t)-U_M(t)\|_\infty=0\,.
\end{equation}
\end{theorem}
The theorem is proved in Section \ref{sec.ab}. There is no restriction on the sign of the solution. By interpolation, we can are easily obtain rates in all $L^q$ spaces, $1<q<\infty$, see for instance examples in \cite{Vaz17}. Of course, for $M=0$ we just say that  $ \|u(t)\|_1$ goes to zero.

It is interesting to interpret Theorem \ref{thm.ab1} in terms of the rescaled variables defined }in Subsection \ref{ssec.ssv}, see formula \eqref{eq.rescflow}. Then we may rephrase the result as saying that $v(y,\tau)$ converges to the equilibrium state $F_M(y)$ of the flow equation \eqref{eq.resc}. In other words, $F_M$ attracts along the rescaled  flow all finite-mass solutions with the same mass.

Section \ref{sec.gharn} settles the question of positivity. It contains the study of two-sided global estimates for nonnegative solutions with compactly supported and bounded data, where the fundamental solution plays a key role, see Theorem \ref{thm.GH}. They are known as {\sl global Harnack inequalities}, though they depend on the information on the initial data.

The last sections contain additional information. Thus, the existence of the source-type solution in a bounded domain is shown in Section~\ref{sec.st-sol}. It is not relevant for the long-time behaviour in the setting of bounded domains, hence it loses interest. Section \ref{sec.limitcases} examines the limit cases $s\to 1$, $s\to 0$,  and $p\to 2$. A final section contains comments on the importance of fundamental solutions in the related literature, followed by other comments and open problems.

\medskip

\noindent {\sc Notations}. We sometimes write a function $u(x,t)$ as $u(t)$ or $u$  when one some of the variables can be safely understood. We use the notation $u_+=\max\{u,0\}$. The letters $\alpha$ and $\beta$ will be fixed at the values given in the self-similar formula \eqref{eq.sse1}. We also use the symbol $\|u\|_q$ as shortened notation for the norm of $u$ in the $L^q$ space over the corresponding domain when no confusion is to be feared. We denote the duality product in $L^q\times L^{q'}$, with $q$ and $q'$ dual exponents, by $\langle \cdot , \cdot\rangle$. For a function $u(x)\ge 0$ we call mass or total mass the integral $\int_{\ren} u(x)\,dx$, either finite or infinite.
For signed functions that integral does not coincide with the $L^1$ norm, so the use of the term is only justified by analogy and usually refers to the $L^1$ norm.  Finally, the sign $\approx$ means equivalence up to a constant positive factor, while $\sim$ means equivalence with limit 1.

\section{Basic theory}\label{sec.basic}

\smallskip

We establish well-posedness of Problem \eqref{frplap.eq}-\eqref{init.data} in different functional spaces, starting by the consideration of the equation as a gradient flow in $L^2(\ren)$. We obtain unique strong solutions that are $C^\delta$-H{\"o}lder continuous and space and time for seme $\delta>0$, and decay as expected by dimensional considerations. We give a detailed account of  the main qualitative and quantitative properties, some of them correspond to known work done in bounded domains, but some are particular to the whole-space setting. Some of the results of the section are new in the literature.

\subsection{Existence and uniqueness}

We can solve the evolution problem for equation \eqref{frplap.eq} with initial data  $u_0\in L^2(\ren)$  by using the fact that the equation is the gradient flow of a maximal monotone operator associated to the  convex functional \eqref{Jsp1}, see for instance \cite{MazRT, Puhst, Vaz16}. Much of the general theory is common to all cases $1<p<\infty$, but since a number of important details differ  for $p<2$ from what is said below, we concentrate on $p>2$. The domain of the operator is
\begin{equation}
D_2({\mathcal L}_{s,p})= \{\phi\in L^2(\ren): \ {\mathcal J}_{s,p}(u)<\infty, \ {\mathcal L}_{s,p}u\in L^2(\ren)\}.
\end{equation}
Well known theory implies that for every initial $u_0\in L^2(\ren)$ there is a unique strong solution $u\in C([0,\infty): L^2(\ren))$, that we may call the semigroup solution. Strong solution means that $u_t$ and $ {\mathcal L}_{s,p}u\in L^2(\ren)$ for every $t>0$, and the equation is satisfied a.e in $x$ for every $t>0$. The semigroup is denoted as $S_t(u_0)=u(t)$, where $u(t)$ is the solution emanating from $u_0$ at time $0$. Typical a priori estimates for gradient flows follow, cf.  \cite{BrBk73, Ko67}.  The next results are part of the standard theory:
\begin{equation}
\frac12\frac{d}{dt}\|u(t)\|_2^2= -\langle  {\mathcal L}_{s,p}u(t),  u(t)\rangle=-p\,
{\mathcal J}(u(t)),
\end{equation}
where ${\mathcal J}={\mathcal J}_{p,s}$ as in the introduction, and also
\begin{equation}
\frac{d}{dt}{\mathcal J}(u(t)) =\langle  {\mathcal L}_{s,p}u(t), u_t(t)\rangle=-\|u_t(t)\|^2_2,
\end{equation}
where integrals and norms are taken in $\ren$. It follows that both $\|u(t)\|_2$ and $J(u(t))$ are decreasing in time, and we get the easy estimate
$$
{\mathcal J}(u(t)\le \|u_0\|_2^2/(2pt)
$$
for every $t>0$. See other properties below.

\noindent $\bullet$  Moreover, for given $p>1$ (the index of the operator) and every $1\le q\le \infty$, the $L^q$ norm of the solution is non-increasing in time.  We can extend the set of solutions to form a continuous semigroup of contractions in $L^q(\ren)$ for every $1\le q< \infty$: for every $u_0\in L^q(\ren)$ there is a unique strong solution such that $u\in C([0,\infty): L^q(\ren))$.  The class of solutions can be called the {\sl $L^q$ semigroup} for equation \eqref{frplap.eq} posed in $\ren$. These $q$-semigroups coincide on their common domain. The Maximum Principle applies, and more precisely $T$-contractivity holds in the sense that for two solutions $u_1, u_2$ and any $q\ge 1$ we have
\begin{equation}
\|(u_1(t)-u_2(t))_+\|_q \le \| (u_1(0)-u_2(0))_+\|_q.
\end{equation}
This implies that we have an ordered semigroup for every $q$ and $p$.  An operator with these properties in all $L^q$ spaces is called completely accretive, cf. \cite{BC91}. We can also obtain the solutions by Implicit Time Discretization, cf. the classical references \cite{CL71, Ev78}. The word mild solutions is used in that context, but mild and strong solutions coincide by uniqueness. The operator is also accretive in $L^\infty$, and this allows to generate a semigroup of contractions in $C_0(\ren)$ the set of continuous functions that go to zero at infinity.

\noindent $\bullet$  This part of the theory can be done for solutions with two signs, but we will often reduce ourselves in the sequel to nonnegative data and solutions. Splitting the data into positive and negative parts most of the estimates apply to signed solutions. To be precise,
for a signed initial function $u_0$ we may consider its positive part, $u_{0,+}$
and its negative part $u_{0,-}= -u_0+ u_{0,+}=-\max\{-u_0, 0\}$. Then, both $u_{0,+}$ and $u_{0,-}$ are nonnegative and $-u_{0,-}\le u_0\le u_{0,+}$. It follows from the comparison property of the $L^q$ semigroups that
$$
-S_t(u_{0,-})\le S_t(u_0)\le S_t(u_{0,+}).
$$
Therefore, we may reduce many of the estimates to the case of nonnegative solutions.

\noindent $\bullet$ An alternative construction approach is to prove that the solutions in $\ren$ are obtained as limits of the solutions of the Dirichlet problem posed in expanding balls $\Omega_R=B_R(0)$, as constructed in \cite{MazRT, Vaz16}. For nonnegative solutions with a common initial datum this limit is monotone in $R$. The proof  that the two ways of construction give the same solutions is easy in the nonnegative case and will be omitted. In this way the $L^q$ semigroups are obtained as limit as the ones on bounded domains and the many properties, like $L^q$ boundedness, contractivity or comparison are inherited.

\subsection{Scaling transformations}\label{ssec.scaling}

In our study we will use the fact that the equation admits a one-parameter scaling group that conserves the mass of the solutions. Thus, if $u$ is a weak or strong solution of the equation, then we obtain a family of solutions of the same type, $u_k={  \mathcal T}_k u$, given by
\begin{equation}\label{scal.trn1}
{\mathcal T}_k u(x,t)= k^{N} u(kx, k^{N(p-2)+sp}t)
\end{equation}
for every $k>0$.  This scaling transformation can be combined with a second one that keeps invariant the space variable
\begin{equation}\label{scal.trn2}
{\widehat{\mathcal  T}}_M u(x,t)= M u(x, M^{p-2}t)
\end{equation}
for every $M>0$. This one can be used to reduce the calculations to solutions with unit mass, $M=1$. Together, these transformations form the two-parameter scaling group under which the equation is invariant.

Let us point out that the set of solutions of the equation is invariant under a number of isometric transformations, like: change of sign: $u(x,t)$ into $-u(x,t)$, rotations and translations in the space variable, and translations in time. They will also be used below.

\subsection{A priori bounds}

\noindent $\bullet$ Our operator is homogeneous of degree $d=p-1>1$ in the sense that ${\mathcal L}_{s,p}(\lambda\,u)=\lambda^{p-1} {\mathcal L}_{s,p}u$. Using the general results by  B{\'e}nilan-Crandall \cite{BC81b}  for suitable homogeneous operators in Banach spaces, we can prove the a priori bound
\begin{equation}
(p-2)t u_t> - u,
\end{equation}
which holds for all nonnegative solutions, in principle in the sense of distributions. This a priori bound is quite universal, independent of the particular nonnegative solution. It is based on  the scaling properties and comparison.  Therefore, we have almost  monotonicity in time if $u\ge0$. In particular, if a strong solution is positive at a certain point $x_0$ at $t=t_0$, then for all later times $u(x_0,t)>0$. This is called {\sl conservation of positivity} (we recall that it holds for nonnegative solutions).

Combined with the decay  of the space integral in time, we conclude another interesting result for nonnegative solutions:
\begin{equation}
\|u_t(\cdot,t)\|_{L^1(\ren)}\le C_p \|u(\cdot,t)\|_{L^1(\ren)}\,t^{-1}.
\end{equation}

\noindent $\bullet$ On the other hand, paper \cite{BC81b} also implies the estimate for all $p>2$ and  $q>1$ we have
\begin{equation}\label{form.brh}
\|u_t\|_q\le \frac 2{(p-2)t}\|u_0\|_q
\end{equation}
for every $1< q\le \infty$. For $q=1$ it is formulated as different quotients.
Formula {\eqref{form.brh} is valid for all signed solutions.

\subsection{Energy estimates}

\noindent $\bullet$ As we have seen before, for solutions with data in $L^2(\ren)$ and for all times $0\le t_1<t_2$ we  have the identity
\begin{equation}\label{en.est}
\int_{\ren} u^2(x,t_1)dx-\int_{\ren} u^2(x,t_2)\,dx=2\int_{t_1}^{t_2}
\int_{{\ren}}\int_{{\ren}} |u(x,t)-u(y,t)|^p\,d\mu(x,y)dt\,,
\end{equation}
\noindent where ${d\mu(x,y)}=|x-y|^{-(N+sp)}dxdy $. In the sequel we omit the domain of integration of most space integrals when it is $\ren$, and the time interval when it can easily understood from the context.

We point out that this estimate shows that solutions with $L^2(\ren)\cap L^p(\ren)$ data belong automatically to the space $L^p(0,\infty: W^{s,p}(\ren))$.

\medskip

\noindent $\bullet$ Arguing in the same way, for solutions with data in $L^q(\ren)$  with $q>1$ and  $0\le t_1<t_2$ we have for nonnegative solutions
\begin{equation}\label{en.est.q}
\begin{array}{c}
\displaystyle \int u^q(x,t_1)dx-\int u^q(x,t_2)\,dx=\\ [6pt]
\displaystyle q\iiint |u(x,t)-u(y,t)|^{p-2}\langle (u(x)-u(y)), (u^{q-1}(x,t)-u^{q-1}(y,t))\rangle \,d\mu(x,y)dt\,,
\end{array}
\end{equation}
with integration in the same sets as before. We use the inequality
\begin{equation}
(a-b)^{p-1}(a^{q-1}-b^{q-1})\ge C(p,q) \,|\,a^{(p+q-2)/p}-b^{(p+q-2)/p}\,|^p
\end{equation}
which is valid for all $a>b>0$ and $p,q>1$. This inequality is also true when $b\ge a>0$ by symmetry, and when $a$ and $b$ have different signs in an elementary way.  We get the new inequality
\begin{equation}\label{en.est.q}
\begin{array}{c}
\displaystyle C(p,q)\int \iint |u(x,t)^{(p+q-2)/p}-u(y,t)^{(p+q-2)/p}|^{p} \,d\mu(x,y)dt\\ [6pt]
\displaystyle \le \int u^q(x,t_1)dx-\int u^q(x,t_2)\,dx,
\end{array}
\end{equation}
which applies the solutions of the $L^q$ semigroup, $q>1$. This gives a precise estimate of the dissipation of the $L^q$ norm along the flow.

Case of signed solutions. The above results hold on the condition that we use the notation $a^{p-1}$ to mean $|a|^{p-2}a$ and so on (this is a usual convention). The equality to prove is
\begin{equation}\label{en.est.q}
\begin{array}{c}
\displaystyle \int |u|^q(x,t_1)\,dx-\int |u|^q(x,t_2)\,dx=\\ [6pt]
\displaystyle q\iiint \langle |u(x)-u(y)|^{p-2}(u(x)-u(y)) , (u^{q-1}(x,t)-u^{q-1}(y,t))\rangle \,d\mu(x,y)dt\,,
\end{array}
\end{equation}
and the dissipation estimate is also true in this case.

Note that these estimates can be obtained as limit of the ones already obtained for the problem posed in a bounded domain.

\subsection{Difference estimates}\label{ss.diff.est}

It is well known that the semigroup is contractive in all $L^q$ norms, $1\le q\le \infty$. At some moments we would like to know how the norms of the difference of two solutions decrease in time. Such decrease is called {\sl dissipation}. We present here the easiest case, decrease in $L^2$ norm.

\noindent {\bf $L^2$ dissipation.} For solutions with data in $L^2(\ren)$ and times $0\le t_1<t_2$ we  have the identity for the difference of two solutions $u=u_1-u_2$:
\begin{equation}\label{diff.est.2}
\begin{array}{c}
\displaystyle\int_{\ren} u^2(x,t_1)\,dx-\int_{\ren} u^2(x,t_2)\,dx \\=
\displaystyle 2\int_{t_1}^{t_2}\int_{{\ren}}\int_{{\ren}}
\displaystyle \left(|u_1(x,t)-u_1(y,t)|^{p-2}(u_1(x,t)-u_1(y,t)) \right. \\[10pt]
 \displaystyle \left.-|u_2(x,t)-u_2(y,t)|^{p-2}(u_2(x,t)-u_2(y,t))\right)\\[8pt]
  \displaystyle (u_1(x,t)-u_2(x,t)-u_1(y,t)+u_2(y,t))\,d\mu(x,y)dt\,,
\end{array}
\end{equation}
\noindent where ${d\mu(x,y)}=|x-y|^{-(N+sp)}dxdy $ as before. Putting $a=u_1(x,t)-u_1(y,t)$ and $b=u_2(x,t)-u_2(y,t) $ and using the numerical inequality as before we bound below the last integral by
$$
C(p) \iiint  |(u_1(x,t)-u_1(y,t))^{p/2}-(u_2(x,t)-u_2(y,t))^{p/2}|^2\,d\mu(x,y)dt\\.
$$
This is an estimate of the $L^2$-dissipation of the difference $u=u_i-u_2$.

Later on, we  will need the expression of the $L^1$ dissipation in the study of the asymptotic behaviour, but we will postpone it until conservation of mass is proved.

\subsection{Boundedness for positive times. Continuity}

\noindent $\bullet$ An important result valid for many nonlinear diffusion problems with homogeneous operators is the so-called $L^1$-$L^\infty$ smoothing effect. In the present case we have

\begin{theorem} For every solution with initial data $u_0\in L^1(\ren)$ we have
 \begin{equation}\label{smoot.effct}
 |u(x,t)|\le C(N,p,s) \|u_0\|_1^{\gamma}\,t^{-\alpha}\,,
 \end{equation}
 with exponents  $\alpha=N\beta$, $\gamma=sp\beta$ and $\beta=1/(N(p-2)+sp$.
  \end{theorem}

The exponents are  given by the scaling rules (dimensional analysis). The result has been recently proved by Bonforte and Salort \cite[Theorem 5.3]{BS2020} where an explicit value for the constant $C(N,p,s)$ is given. It can also be derived as a consequence of the results of   Str{\"o}mqvist \cite{Strom2019}. Note that this formula has to be invariant under the scaling transformations of Subsection \ref{ssec.scaling}. For reference to the similar result in a number of similar nonlinear diffusion theories, including linear and fractional heat equation, porous medium and its fractional versions, $p$-Laplacian, and so on, cf. for instance \cite{Vaz17}.

\noindent $\bullet$ Once we know that solutions are bounded, we can prove further  regularity. We can rely on Theorem 1.2 of \cite{BLS2019} by Brasco-Lindgren-Str{\"o}mqvist that we state in short form as follows:

\begin{theorem} Let $\Omega\subset \ren$  be a bounded and open set, let $I = (t_0, t_1]$, $p \ge 2$ and $0 < s < 1$. Suppose that $u$ is a local weak solution of  \eqref{frplap.eq} in the cylinder  $Q= \Omega\times I$
such that is it bounded in the sense that $u\in  L^\infty_{loc}(I;L^\infty(\ren))$. Then, there exist positive constants $\Theta(s,p)$ and $\Gamma(s,p)$ such that
$$
u\in C^\delta_{x,loc}(Q)\cap C^\gamma_{t,loc}(Q)
$$
for every  $0<\delta<\Theta$ and $0<\gamma <\Gamma$. Moreover, the  H{\"o}lder bounds in both space and time are uniform in any cylinder $Q'=B_R(x_0)\times I'$  strictly included in  $Q$, and they depend only on $N,s,p$, the distance of $Q'$ to the parabolic boundary of $Q$ and on the norm of $u$ in $L^\infty(\ren\times I')$.
 \end{theorem}

Explicit values for $\Theta$ and $\Gamma$ are given in \cite{BLS2019}. We can check that the conditions of this theorem apply to our setting whenever $t_0>0$, hence we have

\begin{corollary}
The solutions of our evolution problem \eqref{frplap.eq}-\eqref{init.data} are uniformly  H{\"o}lder continuous in space with exponent $\delta<\Theta$ and in time with exponent $\gamma <\Gamma$, always for $t\ge t_0>0$.
 \end{corollary}

For completeness we recall a number of previous papers on the elliptic equation ${\mathcal L}_{s,p}u(t)= f$ that proved different results on continuity of solutions of the elliptic version under assumptions on $f$. Let us quote Kuusi-Mingione-Sire \cite{KMSire} who first proved continuity for $sp<N$, Lindgren \cite{Lind16} who proved H{\"o}lder continuity for continuous $f$, Iannizzotto et al. \cite{IannMS2014} who proved H{\"o}lder regularity for bounded $f$ with $u=0$ outside of $\Omega$ and finally Brasco-Lindgren-Schikorra \cite{BLS2018} who proved H{\"o}lder regularity for $f\in L^q_{loc}$ with $q>N/sp$, $q\ge 1$. This last result was the basis of an alternative but more complicated former proof we had for our corollary.

\subsection{Positivity of nonnegative solutions}\label{ssec.pos}

Nonnegative strong solutions of equation \eqref{frplap.eq} enjoy the property of strict positivity at least in the almost everywhere sense. Indeed, at every point $(x_0,t_0)$ where  a  solution reaches the minimum value $u=0$ and
$ ({\mathcal L}_{s,p}u)(x_0,t_0)$ exists, then it must be strictly negative according to the formula for the operator. On the other hand, if $u_t$ exists it must to zero. From this contradiction we conclude that a.e. $u(x,t)$ must be positive. By the already proved conservation of positivity, for any $t>t_0$ we have  $u(x,t_0)>0$ for a.\,e. $x\in \ren$.

Since we already know that the nonnegative solution is continuous, then $u$ is positive everywhere unless if it is zero everywhere. Quantitative positivity will be discussed in Section \ref{sec.pos.tail}.

\subsection{Comparison via symmetries. Almost radiality}

The Aleksandrov symmetry principle \cite{Alek} has found wide application in elliptic and parabolic linear and nonlinear problems. An explanation of its use for the Porous Medium Equation is given in  \cite{Vazpme07}, pages 209--211. In the parabolic case it says that whenever an initial datum can be compared with its reflection with respect to a space hyperplane, say $\Pi$, so that they are ordered, and the equation is invariant under symmetries, then the same space comparison applies to the solution at any positive time $t>0$.

The result has been applied to elliptic and parabolic equations of Porous Medium Type involving the fractional Laplacian in \cite{VazBar2014}, section 15. The argument of that reference can be applied in the present setting. We leave the verification to the reader-

The standard consequence we want to derive is the following

\begin{proposition} \label{Alek}Solutions of our Cauchy Problem having compactly supported data in a ball $B_R(0)$ are radially decreasing  in space for all $|x|\ge 2R$. Moreover, whenever $|x|>2R$ and $|x'|<|x|-2R$, then we have $u(x,t)\le u(x',t)$
for all $t>0$.\end{proposition}

\subsection{On the fundamental solutions}\label{ssec.bddness}

The existence and properties of the fundamental solution of Problem \eqref{frplap.eq}-\eqref{init.data} are a main concern of this paper. We expect it  to be unique, positive and self-similar for any given mass $M>0$. Self-similar solutions have the form
\begin{equation*}\label{eq.ssf1}
U(x,t;M)=t^{-\alpha}F(x\,t^{-\beta};M)
\end{equation*}
(more precisely, this is called direct self-similarity). Substituting this formula into equation  \eqref{frplap.eq}, we see that time is eliminated as a factor in the resulting formula on the condition that: $\alpha+1=(p-1)\alpha+\beta sp$. We also want integrable solutions that will enjoy the mass conservation property, which implies $\alpha=N\beta$. Imposing both conditions, we get
\begin{equation*}\label{eq.sse}
\alpha=\frac{N}{N(p-2)+sp}, \quad \beta=\frac{\alpha}{N}=\frac{1}{N(p-2)+sp}\,,
\end{equation*}
as announced in the Introduction. Note that for $p>2$ we have $\alpha> N/2s$ and $\beta>1/2s$. The profile function $F(y;M)$ must satisfy the nonlinear stationary fractional equation
\begin{equation}\label{fsprof}
{\mathcal L}_{s,p} F= \beta \,\nabla\cdot (yF)\,.
\end{equation}
Cf. a similar computation for the Porous Medium Equation in \cite{Vazpme07}, page 63. Using rescaling ${\widehat {\mathcal  T}}_M $, we can reduce the calculation of the profile to mass 1 by the formula
$$
F(y;M)=M^{sp\beta}F(M^{-(p-2)\beta}y;1).
$$
In view of past experience with $p=2$, we will look for  $F$ to be radially symmetric, monotone nonincreasing in $r=|y|$, and positive everywhere with a certain behaviour as $|y|\to \infty$.

We have proved that all solutions with $L^1$ data at one time will be uniformly bounded and continuous later on. Thus,  $F$ must  be bounded and continuous. Moreover, bounded solutions have a bounded $u_t$ for all later times. In the case of the fundamental solution,  this means that $r^{1-N}(r^NF(r))'$ is bounded, hence $rF'$ is bounded,  and $F$ is regular for all $r>0$.

The self-similar fundamental solution must take a Dirac mass as initial data, at least in the sense of initial trace, i.\,e., $u(x,t)\to M\delta(x)$ as $t\to 0$ in a weak sense. It  will be invariant under the scaling group  ${{\mathcal  T}}_k $ of Subsection \ref{ssec.scaling}. All of this will be proved in the sequel. The detailed statement is contained in Theorems \ref{thm.exfs} and \ref{thm.exfs2} and whole proofs follow there.


\subsection{Self-similar variables}\label{ssec.ssv}
In several instances in the sequel it will be quite convenient to pass to self-similar variables, by zooming the original solution according to the self-similar exponents \eqref{eq.sse}.
More precisely, the change is done  by the formulas
\begin{equation}\label{eq.rescflow}
u(x,t)=(t+a)^{-N\beta}v(y,\tau) \quad y=x\,(t+a)^{-\beta}, \quad \tau=\log(t+a)\\,
\end{equation}
with $\beta=(N(p-2)+sp)^{-1}$, and any $a>0$. It implies that  $v(y,\tau)$ is a solution of the corresponding PDE:
\begin{equation}\label{eq.resc}
\partial_\tau v + {\mathcal L}_{s,p} v -\beta \nabla\cdot(y\, v)=0\,.
\end{equation}
This transformation is usually called \sl continuous-in-time rescaling \rm to mark the difference with the transformation with fixed parameter \eqref{scal.trn1}.

Note that the rescaled equation does not change with the time-shift $a$ but the initial value in the new time does,
$\tau_0=\log(a)$, If $a=0$ then $\tau_0=-\infty$ and the $v$ equation is defined for $\tau\in \re$. The mass of the $v$ solution at new time $\tau\ge \tau_0$ equals that of the $u$ at the corresponding time $t\ge 0$.

Sometimes $\tau$ is defined as $\tau=\log((t+a)/a)$ without change in the equation. It is just a displacement in the new time, but it is important to take it into account in detailed computations.

Denomination: for convenience we will sometimes refer below to the solutions of the rescaled equation \eqref{eq.resc} as $v$-solutions, while the original ones are called $u$-solutions. We hope this is a minor licence.

\medskip

\section{Barrier construction and tail behaviour}\label{sec.barr}

Here we will construct an upper barrier $\widehat u(x,t)$ for the solutions  of the  Cauchy problem with suitable initial data. The barrier will be needed later in the construction of the fundamental solution as limit of approximations with the same mass as the initial Dirac delta.  We will only need to consider nonnegative data and solutions. Besides, it will be enough to do it for bounded radial functions with compact support as initial data, and then use some  comparison argument  to  eliminate the restrictions of radial symmetry and compact support. The barrier will be radially symmetric, decreasing in $|x|$ and will have  behaviour $\widehat u(x,t)=O(|x|^{-N-sp})$ for very large $|x|$.  Note that such behaviour is integrable at infinity.

\noindent {\bf Rescaling.} We will work with the rescaled solution  and the equation \eqref{eq.resc} introduced in Subsection \ref{ssec.ssv}. Translating previous a priori bounds for the original equation into the present rescaled version, we see that all the rescaled solutions are bounded $v(r,t) \le A((t+a)/t)^{N\beta}$. As a consequence of finite mass, radially symmetry and monotonicity in the radial variable, we also get a bound of the form
$$
v(|y|,\tau)\le B \,|y|^{-N}\,,
$$
 and the decay at infinity  is uniform in time, $B$ does not depend on time.  $A$ and $B$ depend on the mass of the solution.

Therefore, we only need to refine the latter estimate for large $r$ so that we get an integrable barrier in an outer region $r\ge R_1\gg 0$.
 We use the notation $r=|y|>0$ in this section where we work with self-similar variables.

\noindent {\bf Construction.} The upper barrier  we consider in self-similar variables   will be stationary in time, $\widehat v(y)=G(r)$, $r=|y|$.  The barrier will have the form of an inverse power in the far field region. We need to compare $v(r,t)$ with $G r)$ in an outer domain, and make a correction of the solution in the near field so that the needed comparison works.

 To be precise, the barrier will be defined by different expressions in three regions: We select two radii $1<R<R_1$. We take $R_1>R$, in fact much larger than $R$. For  $r>R_1$ \ we take the form
\begin{equation}\label{barrier}
G(y)=C_1 r^{-(N+\gamma)}, \quad r=|y|\,.
\end{equation}
We need a  $\gamma>0$, we will later make the choice $\gamma=sp$ that will turn out to be sharp. For $r\le R$ it is smooth and proportional to $A$. But constant equal to $A$ is fine. Finally, in the intermediate region $R<r<R_1$
we put
\begin{equation}
G(r)\sim C_2 r^{-N}.
\end{equation}
We have to glue these regions: $A= C_2R^{-N}$, $C_2R_1^{-N}= C_1 R_1^{-N-sp}$. We can do it in a smooth way, but it is not necessary the details are not important.
\medskip

\noindent {\bf Supersolution in the outer region.} The main difficulty lies in the comparison in the domain that is the exterior of a big ball, and for a long interval of time, \ $Q=\{r> 2R_1\}\times (0,T)$. We want to prove that given a solution $v$ with small initial data, then \ $v\le \widehat v$ by in $Q$ by using the equation in rescaled form plus the interior and initial conditions.
The most difficult part is to prove the supersolution condition for the equation in $Q$.

\begin{lemma} \label{barr.lem} If \ ${\widehat v}(y)=G(r)$ is defined as above, and the positive constants $A,$ $R,$ $R_1,$ $C_1$, $C_2$, $C_3$ are suitably chosen, see below, then
 for all $r\ge 2R_1$ we have
\begin{equation}
{\mathcal L}_{s,p} {\widehat v}-\beta \, r^{1-N}(r^N{\widehat v})_r\ge 0\,.
\end{equation}
\end{lemma}

\noindent {\sl Proof.} The computation of the right-hand side is immediate:
\begin{equation}
-\beta \, r^{1-N}(r^N{\widehat v} )_r=\beta\gamma C_1r^{-N-\gamma}>0  \,,
\end{equation}
which has a good sign. This quantity must control all possible negative terms.

 On the other hand, ${\mathcal L}_{s,p} \widehat v$  may be negative. We have a series of partial estimates of the  contribution of different regions that are to be compared against  the previous bound.

(i) The first term to comes from the influence of the inner core $\{r<R\}$ where we have  ${\widehat v} \approx A$. We get for the contribution from this region to the integral ${\mathcal L}_{s,p} \widehat v$  the quantity
$$
I_1:=
\left.{\mathcal L}_{s,p} \widehat  v\,\right]_1\sim -A^{p-1}R^N r^{-N-sp}
$$
For the moment we need $A^{p-1}R^N r^{-N-sp}\le C_1r^{-N-\gamma}$, that holds if $\gamma\ge sp$ and $A^{p-1}R^N \le C_1$. We fix $\gamma=sp$  in the sequel. We need $A^{p-1}R^N\le \ve_1\, C_1$.

(ii) We still need to calculate the contribution of the remaining regions. Let us fix the point $r=r_0>2R_1$. The contribution of the region $\{r>r_0\}$ need not be counted if we do not like to since the integrand is positive, see the formula. In the ball $D_2=B_0= B_{r_0}(r_0)$ that is not strictly contained in the annulus,  we have
$$
I_2:=\left.-{\mathcal L}_{s,p} \widehat  v(r_0)\,\right]_2 =
\int_{B_0}\frac{\Phi_p(\widehat v(x))-\Phi_p(\widehat v(x+z))}{|z|^{N+sp}}\,dz
$$
where $z=y-y_0$, $|y_0|=r_0$. Since $\widehat  v$ is a $C^2$ function without critical points in $B_0$ the integral converges by the results of \cite{KKL}, Section 3. This  takes into account the cancellations of differences at points located symmetrically w.r.to $y_0$. A way of doing this here is to use the equivalent symmetrized form
$$
I_2=\frac12  \int_{B_0}\frac{2\Phi_p(\widehat v(x))-\Phi_p(\widehat v(x+z)-\Phi_p(\widehat v(x-z))}{|z|^{N+sp}}\,dz
$$
that we can write as
$$
I_2=\frac{C_1^{p-1}}2 \int_{B_0}\frac{2f(x))-f(x+z)-f(x-z))}{|z|^{N+sp}}\,dz
$$
with $f(x):=|x|^{-(p-1)(N+sp)}$. We have
$$
f'(r)= c' \, r^{-(N+sp)(p-1)-1},  \quad \widehat  v''(r)= c''\, r^{-(N+sp)(p-1)-2}
$$
and we have the following estimate in $B_0$:
$$
|2f(x))-f(x+z)-f(x-z))|\le c |D^2f(x')||z|^2,
$$
where $x'$ is an intermediate point. We thus get with $\rho=|z|=|y-y_0|$
$$
|I_2|\sim c  C_1^{p-1} \int_{0}^{r_0/2}  r_0^{-(N+sp)(p-1)-2}\rho^{1-sp}\,d\rho
$$
Since $sp<2$ we finally arrive at the estimate
$$
|I_2|\le c\,C_1^{p-1}\,r_0^{-\gamma'}, \quad
\gamma'= (N+sp)(p-1)+sp.
$$
Now,
$$
\gamma'-(N+sp)= (N+sp)(p-2)+sp>0
$$
Therefore, the calculation enters our scheme if $C_1^{p-2}\le \ve_2 r_0^{\gamma'-N-sp}$.
\nc

(iii) We have to calculate the contribution of the rest of the annulus $\{R<r<r_0\}$.  In the region $D_3=\{R_1<r<r_0\}\setminus B_{r_0}(r_0)$, we get the contribution:
$$
I_3:=
\left.-{\mathcal L}_{s,p} \widehat  v(r_0)\,\right]_3 \le\int_{R_1}^{r_0/2} c  C_1^{p-1} r^{-(N+sp)(p-1)}r_0^{-(N+sp)}r^{N-1}dr
\sim c\,C_1^{p-1}\,R_1^{-\gamma_1} r_0^{-(N+sp)},
$$
with
$$
\gamma_1= (N+sp)(p-1)-N=N(p-2)+ps(p-1)>0.
$$
We need $C_1^{p-2}\le \ve_3 R_1^{\gamma_1}$.

(iv) Finally, for $D_4=\{R<r<R_1\}$
$$
I_4:=
\left.-{\mathcal L}_{s,p} \widehat  v(r_0)\right|_4 \le\int_{R}^{R_1} c  C_2^{p-1} r^{-N(p-1)}r_0^{-(N+sp)}r^{N-1}dr
\sim c\,C_2^{p-1}\,R^{-N(p-2)} r_0^{-(N+sp)},
$$
so that we need $C_2^{p-1}\le \ve_3 C_1 R^{N(p-2).}$

\noindent $\bullet $ List of inequalities
$$
A^{p-1}R^N\le \ve_1\, C_1, \quad
C_2^{p-1}\le \ve_4 C_1 R^{N(p-2),}\quad
C_1^{p-2}\le \ve_2 R^N R_1^{\gamma'-N-sp}, \quad
C_1^{p-2}\le \ve_3 R_1^{\gamma_1},
$$
together with $A= C_2R^{-N}$, $C_2R_1^{sp}= C_1 $. Using the equalities, the former list becomes the series of conditions
\begin{equation}
\begin{split}
&A^{p-2}\le \ve_1 R_1^{sp}, \quad
C_2^{p-2}\le \ve_4  R^{N(p-2)}\,R_1^{sp} \\
&C_2^{p-2}\le \ve_2 R^N R_1^{\gamma'-N+sp(p-3)}, \quad
C_2^{p-2}\le \ve_3 R_1^{\gamma_1+sp(p-2)},
\end{split}
\end{equation}

After choosing $A$ and $C_2$ we fix $R$ with $A= C_2R^{-N}$. Then we need a large $R_1$  to satisfy the rest of the inequalities. The construction is done and becomes a supersolution of the equation in that region. By usual comparison for fractional equations, it will be on top in this region of any solution with conveniently small initial data with compact support. But that question will be discussed next at the global level. \qed

\medskip

\begin{theorem}[Global barrier]\label{thm.barr} Let us define a function $G(|y|)$ by the previous recipe, with constants as large as prescribed. Let $v_0\ge 0$, is integrable with mass 1, $v_0(y)$ is radially symmetric and decreasing w.r.to the radius. Then, if we assume that $v_0\le G$ in $\ren$, it follows that  $v(\tau)\le G$ for all times.
\end{theorem}

\noindent {\bf Proof.} We point out that assuming $\int u_0(x)\,dx=M=1$ is no restriction because of the mass changing transformation. We may also assume that $u_0$ is more regular to justify the comparison calculations, by density this is no problem. We define the barrier with a different formula in every region as prescribed above. We have to check that $v(y,\tau)$ may not touch $G$ from below in any of them

(i) In the intermediate region  $\mathcal{D}_2=\{y: R\le |y|\le 2R_1\}$ we have $G(y)=\min\{C_2|y|^{-N}, C_1|y|^{-N-sp}\}$. By the assumption of monotonicity in $|y|$ and the mass assumption we get $v(y,\tau)\le c(N)M\,|y|^{-N}$, so that we have $v<G$ everywhere in $Q_2=\mathcal{D}\times (0,\infty)$ if we take $C_2$ large enough so that $C_2>c(N)$.

(ii) In the inner region $\mathcal{D}_1=\{y:  |y|\le R\}$ we put $G(y)=A$ with $A$ large enough. Let us see that the supersolution condition holds everywhere for the equation :
$$
{\mathcal L}_{s,p} G(y)-\beta \, r^{1-N}(r^N G)_r\ge 0.
$$
Let us make it happen: On one hand, $\beta \, r^{1-N}(r^N G)_r= \beta N G= \alpha A$. On the other hand, $G$ attains it maximum, so it is positive. Moreover,
$$
{\mathcal L}_{s,p} G(y)\ge \int \frac{(G(y)-G(y'))^{p-1}}{|y-y'|^{N+sp}}\,dy' \ge\int_{|y|>R_1|}\frac{(A-G(R_1))^{p-1}}{|y-y'|^{N+sp}}\,dy\ge c(N,c,p)A^{p-1}.
$$
We conclude that ${\mathcal L}_{s,p} G(y)-\beta \, r^{1-N}(r^N G)_r\ge 0$ if $A^{p-2}$ is large enough depending only on $N,s,p$.

(iii) In the outer region $\mathcal{D}_3=\{y: |y|\ge 2R_1\}$ we consider the value  $G(y)=C_1 |y|^{-(N+2s)}$. We need to prove that $G$ is a supersolution of the equation everywhere in the region for all time. This is the difficult part that we have separated as Lemma \ref{barr.lem}.

(iv) The proof of comparison is done in the usual way, we assume that $u_0\in L^2(\ren)$ and $u_0\le G$ and compute
 $$
 \frac{d}{dt}\int (v-G)_+^2\,\dx=\int_{\ren} \left[\beta \, r^{1-N}(r^N (v-G)_r -({\mathcal L}_{s,p} v(y)-{\mathcal L}_{s,p} G(y))\right]\, (v-G)_+\,dx
 $$
Since $(v-G)_+=0$ in $\mathcal{D}_2$ for all times, we can reduce the domain of integration to  $\mathcal{D}_1\cup \mathcal{D}_3$ where the first factor of the integrand is negative. Therefore, we conclude that
$$
 \frac{d}{dt}\int (v-G)_+^2\,\dx\le 0,
 $$
 hence, $v\le G$ for all times if it holds at zero. Comparison can also be done by the viscosity method that is quite intuitive. \qed

 The global barrier can be used to find a rate of decay in space of the solutions which is uniform for bounded mass and some initial  decay at infinity.  The main result is the following

\medskip

\begin{corollary}\label{cor.decay} Let $u$ be a solution with nonnegative, bounded and compactly supported data $u_0$. Then, for every $x\in\ren$, $t>0$ we have
\begin{equation}
u(x,t)\le U(x,t):=(t+1)^{-\alpha} G(|x|\,(t+1)^{-\beta})\,,
\end{equation}
where $G$ is a positive and bounded function such that \ $G(r)\le Cr^{-(N+sp)}$. $C$ depends only on $s,p,N$ and the bounds on the data.
    \end{corollary}

 \noindent {\bf Remark.}   We have taken $a=1$ for convenience since then $\tau_0=0$, $x=y$  and $v(y,0)=u_0(x)$.
 The same formula holds with $(t+a)$ instead of  $(t+1)$ but then $C$ changes.

\section{Mass conservation}\label{sec.masscon}

We now proceed with the mass analysis. The main result is the conservation of the total mass for the Cauchy problem posed in the whole space with nonnegative data.

\begin{theorem}\label{thm.mc} Let $u(x,t)$ be the semigroup solution of Problem \eqref{frplap.eq}, \eqref{init.data}, with $u_0\in L^1(\ren)$, $u_0\ge 0$. Then for every $t>0$ we have
\begin{equation}
\int_{\ren} u(x,t)\,dx=\int_{\ren} u_0(x)\,dx.
\end{equation}
\end{theorem}

Before we proceed with the proof we make two reductions: i)  We may always assume that $u_0\in L^1(\ren)\cap L^\infty(\ren)$ and compactly supported. If mass conservation is proved under these assumptions then it follows for all data $u_0\in L^1(\ren)$ by the contraction semigroup.

We recall that the $L^1$ mass is not conserved in the case of the Cauchy-Dirichlet problem posed in a bounded domain since mass flows out at the boundary. Indeed, the mass  decays in time according to a power rule. On the other hand, mass conservation holds for the most typical linear and nonlinear diffusion problems posed in $\ren$, like the Heat Equation, the Porous Medium Equation or the evolution $p$-Laplacian equation. It also holds for Neumann Problems with zero boundary data posed on bounded domains.

The proof of the theorem is divided into several cases in order to graduate the difficulties.

\subsection{First case: $N<sp$.}\label{ssec.4.1} Here the mass calculation is quite straightforward. We do
a direct calculation for the tested mass. Taking a smooth and compactly supported test function $\varphi(x)\ge0 $, we have for $t_2>t_1>0$
\begin{equation}\label{mass.calc}
\left\{\begin{array}{l}
\displaystyle \left|\int u(t_1)\varphi\,dx-u(t_2)\varphi\,dx\right|\le \iiint
\left|\frac{\Phi(u(y,t)-u(x,t))(\varphi(y)-\varphi(x)}{|x-y|^{N+sp}}\right|\,dydxdt\\[10pt]
\le \displaystyle \left(\iiint  |u(y,t)-u(x,t)|^p\,d\mu(x,y)dt\right)^{\frac{p-1}{p}}
\left(\iiint |\varphi(y)-\varphi(x)|^p\,d\mu(x,y)dt\right)^{\frac{1}{p}}\,,
\end{array}\right.
\end{equation}
with space integrals over $\ren$. Use now the sequence of test functions $\varphi_n(x)=\varphi(x/n)$ where $\varphi(x)$ is a cutoff function which equals 1 for $|x|\le 2$ and zero for $|x|\ge 3$. Then,
$$
\int \int |\varphi_n(y)-\varphi_n(x)|^p\,d\mu(x,y)=
n^{N-sp}\int \int |\varphi(y)-\varphi(x)|^p\,d\mu(x,y)
$$
and this tends to zero as $n\to\infty$. Using \eqref{en.est} we conclude that the triple integral involving $u$ is also bounded in terms
of $\|u(\cdot,t_1)\|_2^2$, which is bounded independently of $t_1$. Therefore, taking the limit as $n\to\infty$ so that $\varphi_n(x)\to 1$ everywhere, we get
$$
\int u(x,t_1)\,dx=\int u(x,t_2)\,dx,
$$
hence the mass is conserved for all positive times for data in $L^2\cap L^1$. The statement of the theorem needs to let $t_1\to 0$, but this can be done thanks to the continuity of solution of the $L^1$ semigroup as a curve in $L^1(\ren)$.

The limit case $N=sp$ also works by revising the integrals, but we get no rate.

\subsection{Case $N\ge  sp$.}\label{ssec.4.2} In order to obtain the mass conservation in this case  we need to use a uniform estimate of the  decrease of the solutions in space so that they help in estimating the convergence of the integral. This will be done by using the barrier estimate that we have proved.

\noindent {$\bullet$} We go back to the first line of \eqref{mass.calc}. The proof relies on some calculations with the multiple integral in that line. We also have to consider different regions. We first deal with exterior region $A_n=\{(x,y): |x|,|y|\ge n\}$,  where recalling \eqref{mass.calc} we have
\begin{equation*}\label{mass.calc2}
\begin{array}{c}
 \displaystyle I(A_n):= \int_{t_1}^{t_2}\iint_{A_n} \frac{|\Phi(u(y,t)-u(x,t))|\,|\varphi_n(y)-\varphi_n(x)|}{|x-y|^{N+sp}}\,dydx\,dt\\[10pt]
\le \displaystyle \left(\iiint  |u(y,t)-u(x,t)|^{p}\,d\mu(x,y)\,dt \right)^{\frac{p-1}{p}}
\left(\iiint |\varphi_n(y)-\varphi_n(x)|^{p}\,d\mu(x,y)dt \right)^{\frac{1}{p}}
\end{array}\end{equation*}
which we write as $I=I_1. I_2$. In the rest of the calculation we omit the reference to the limits that is hopefully understood.

We already know that $I_2 \le C_p\,n^{(N-sp)/p}(t_2-t_1)$. On the other hand,
we want to compare $I_1$ with the dissipation $D_\ve$ of the $L^r$ norm, for $r=1+\ve$. We recall that
\begin{equation*}\label{mass.calc3}
\displaystyle D_\ve= \iiint |(u(y,t)-u(x,t)|^{p-1}\,|u^\ve(y,t)-u^\ve(x,t)|  \,d\mu dt \le C(\ve,p)\int |u|^{1+\ve}(x,t_1),dx\le  C(\ve,p,u_0).
\end{equation*}
Next, we use the elementary equivalence: for all $\ve\in (0,1)$ and all $a,b>0$ we have
$$
|a^\ve-b^\ve| \ge C(\ve) (a-b) (a+b)^{\ve-1}\,.
$$
It follows that
$$
 \displaystyle D_\ve \ge C_\ve \iiint |(u(y,t)-u(x,t)|^{p}\,(|u|^\ve(y,t)+|u|^\ve(x,t))^{\ve-1}  \,d\mu dt\,.
$$
After comparing the formulas, we conclude that
$$
I_1^{p/(p-1)}\le C\,D_\ve\|2u\|_\infty^{1-\ve}\,,
$$
In view of the value of $u$ in the region $A_n$, $u\approx n^{N+sp}$, we have $I(A_n)\le C n^{-\sigma}$ with
$$
\sigma =\frac1{p}\left((N+sp)(p-1)(1-\ve)-(N-sp)\right).
$$
Since $p\sigma = N(p-2)+ sp^2-\ve (N+sp)(p-1)>0$ for $\ve$ small, this gives the vanishing in the limit  $n\to\infty$ of this term that contributes to the conservation of mass. Note that the argument holds for all $p\ge 2$ and $0<s<1$.

\noindent $\bullet$  We still have to make the analysis in the other regions. In the inner region $B_n=\{(x,y): |x|,|y|\le 2n\}$ we get
$\varphi_n(x)-\varphi_n(y)=0$, hence the contribution to the integral \eqref{mass.calc} is zero. It remains to consider the cross regions
$C_n=\{(x,y): |x|\ge 2n ,|y|\le n\}$ and $D_n=\{(x,y): |x|\le n ,|y|\ge 2n\}$. Both are similar so we will look only at $C_n$. The idea is that we have an extra estimate: $|x-y|>n$ so that
\begin{equation*}
 \begin{array}{c}
 \displaystyle I(C_n)\le n^{-(N+sp)}\int_{t_1}^{t_2}\iint_{C_n}|\Phi(u(y,t)-u(x,t)|\,|\varphi_n(y)-\varphi_n(x)|\,dydx\,dt\\[10pt]
\displaystyle  \le Cn^{-(N+sp)}(t_1-t_2)\int dy \int |u(x,t)|^{p-1}\,dy\le Cn^{-(N+sp)}(t_1-t_2)n^N\|u_0\|_{p-1}^{p-1}\,,
\end{array}
\end{equation*}
which tends to zero as $n\to \infty$ with rate $O(n^{-sp})$. Same  for $I(D_n)$. This concludes the proof. Note that these regions overlap but that is no problem. \qed

\medskip

{\bf Signed data.} Theorem \ref{thm.mc} holds also for signed data and solutions. However, the denomination mass for the integral over $\ren$ is physically justified only when $u\ge 0$. For signed solutions the theorem talks about conservation of the whole space integral. The above proof has be  reviewed. Subsection \ref{ssec.4.1} needs no change. As for Subsection \ref{ssec.4.2}, the elementary equivalence has to be written for all $a,b\in\re$
$$
|a^\ve-b^\ve| \ge C(\ve) |a-b| (|a|+|b|)^{\ve-1}\,.
$$

\normalcolor

\subsection{A quantitative positivity lemma}

As a consequence of mass conservation and the existence of the upper barrier we obtain a positivity lemma for certain solutions of the equation.

\begin{lemma}\label{lem.pos}
Let $v$ be the solution of equation \eqref{eq.resc} with initial data $v_0$ such that: $v_0$ is a bounded, nonnegative function with support in the ball of radius $R$, it is radial and radially decreasing, and  $\int v_0(y)\,dy=M>0$. Then
there is a continuous nonnegative function $\zeta(y)$, positive in a ball of radius $r>0$, such
that for every $\tau>0$
\begin{equation}
v(y,\tau)\ge  \zeta(y) \quad \mbox{ for all } \ y\in \ren, \ \tau>0.
\end{equation}
In particular, we may take $\zeta(y)\ge c_1>0$ in $ B_{r_0}(0)$ for suitable $r_0$ and $c_1>0$, to be computed  below.
\end{lemma}

\noindent {\sl Proof.} We know that for every $\tau>0$ the solution $v(\cdot,t)$ will be nonnegative, radial, radially nonincreasing. By Section \ref{sec.barr} there is an upper barrier $G(y)$ on top of $v(y,\tau)$ for every $\tau$. Since $G$ is integrable, for every $\varepsilon>0$ small there is $R(\ve)>0$
such that
$$
\int_{\{|y|>R(\ve)\}} v(y,\tau)\,\dy\le
\int_{\{|y|>R(\ve)\}} G(y\,\dy\le \ve
$$
for all $\tau>0$. Moreover, there is a radius $r_0>0$ such that
$$
\int_{\{|y|<r_0\}} v(y,\tau)\,\dy\le
\int_{\{|y|<<r_0\}} G(y\,\dy\le M/3
$$
for all $\tau>0$. Therefore,
$$
\int_{\{r_0\le |y|\le R(\ve)\}} v(y,\tau)\,\dy\ge M-\ve-M/3 >M/2.
$$
Since $v$ is monotone in $r=|y|$ we have the result
$$
v(r_0,\tau)(R(\ve)^N-r_0^N) \ge c(n)M/2,
$$
hence $v(r,\tau)\ge c_1$ for all $r\le r_0$ and $\tau>0$, with $c_1=c(N,s,p,M,R)$.
Note that the qualitative argument does not depend on the initial $M$ and $R$. \qed


\section{$L^1$ dissipation for differences}\label{sec.diff}
In subsequent sections we will need the very interesting case of the dissipation of the difference $u=u_1-u_2$ in the framework of the $L^1$ semigroup.
We multiply the equation by $\phi=s_+(u_1-u_2)$, where $s_+$ denotes the sign-plus or Heaviside function, and then integrate in space and time. We get in the usual way, with $u=u_1-u_2$, $u_+=\max\{u,0\}$,
\begin{equation}\label{diff.est.1}
\begin{array}{c}
\displaystyle \int u_+(x,t_1)\,dx-\int u_+(x,t_2)\,dx=\int_{t_1}^{t_2}\int s_+(u)u_t\,dx\\ [6pt]
\displaystyle  =\int_{t_1}^{t_2}dt\int ({\mathcal L}_{s,p}  u_1-{\mathcal L}_{s,p}  u_2)\,s_+(u_1-u_2)\,dx =\\ [6pt]
\displaystyle \int_{t_1}^{t_2}dt \iint \left(|u_1(x,t)-u_1(y,t)|^{p-2}(u_1(x,t)-u_1(y,t)) \right. \\[10pt]
 \displaystyle \left.-|u_2(x,t)-u_2(y,t)|^{p-2}(u_2(x,t)-u_2(y,t))\right)\\[8pt]
\left( s_+(u(x,t))-s_+(u(y,t))\right)\,d\mu(x,y)\,.
\end{array}
\end{equation}
We recall that $s_+(u(x,t))=1$ only when $u_1(x,t)>u_2(x,t)$, and
$s_+(u(y,t))=0$ only when $u_1(y,t)<u_2(y,t)$. If we call the last factor in the above display
$$
I=s_+(u_1(x,t)-u_2(x,t))-s_+(u_1(y,t)+u_2(y,t))\,,
$$
we see that $I=1$ if $u_1(x,t)>u_2(x,t)$ and $u_1(y,t)\le u_2(y,t)$. Therefore, on that set
$$
u_1(x,t)-u_1(y,t)>u_2(x,t)-u_2(y,t).
$$
In that case we examine the other factor,
$$
F= |u_1(x,t)-u_1(y,t)|^{p-2}(u_1(x,t)-u_1(y,t)) -|u_2(x,t)-u_2(y,t)|^{p-2}(u_2(x,t)-u_2(y,t))\,,
$$
and conclude that it is positive. The whole right-hand integrand is positive.

In the same way, $I=-1$ if $s_+(u(x,t))=0$ and $s_+(u(y,t)=1$ i.e., only when $u_1(x,t)\le u_2(x,t)$ and $u_1(y,t)>u_2(y,t)$. Then,
$u_1(x,t)-u_1(y,t)<u_2(x,t)-u_2(y,t)$ and $F<0$. The whole right-hand integrand is again positive. We conclude that

\begin{proposition} In the above situation we have the following dissipation estimate:
\begin{equation}\label{diff.est.1b}
\begin{array}{c}
\displaystyle \int (u_1-u_2)_+(x,t_1)\,dx-\int (u_1-u_2)_+(x,t_2)\,dx\\ [8pt]
\ge \displaystyle\iiint_D \left| |u_1(x:y,t)|^{p-2}u_1(x:y,t) -|u_2(x:y,t)|^{p-2}u_2(x:y,t)\right|\,d\mu \,dt.
\end{array}
\end{equation}
where $u_1(x:y,t)=u_1(x,t)-u_1(y,t)$, $u_2(x:y,t)=u_2(x,t)-u_2(y,t)$, and $D\subset \re^2$ is the domain where
$$
\{ u(x,t)>0, \  u(y,t)\le 0\}\cup\{  u(x,t)\le 0, \  u(y,t)>0\}\\,
$$
that includes the whole domain where $u(x,t)\,u(y,t)<0$. There is no dissipation on the set where $u(x,t)\,u(y,t)>0$.
\end{proposition}


\section{Existence of a fundamental solution}\label{sec.fs}

This section deals only with nonnegative solutions unless mention to the contrary. This the first main result.

\begin{theorem}\label{thm.exfs} For any value of the mass $M>0$ there exists a fundamental solution of Problem \eqref{frplap.eq}-\eqref{init.data} having the following properties: (i) it is a nonnegative strong solution of the equation in all $L^q$ spaces, $q\ge 1$, for $t\ge t_0>0$. (ii) It is  radially symmetric and decreasing in the space variable. (ii) It decays in space as predicted by the barrier, $u(t)=O(|x|^{-N+sp})$. (iii) It decays in time $O(t^{-\alpha})$ uniformly in $x$.
\end{theorem}

\noindent {\sc Proof}. We will use the  rescaling method to construct the fundamental solution as a consequence of some asymptotic behaviour as $t\to \infty$. This method has been used in typical nonlinear diffusion problems like the Porous Medium Equation, see \cite{Vazpme07}, and relies on suitable a priori estimates, that are available after the previous sections. The version of the method we use here is the continuous rescaling, that can be of independent interest for the reader.

\noindent $\bullet$  We take an initial datum $\phi(x)\ge 0$ that is bounded, radially symmetric and supported in the ball of radius 1 and has total mass $M=1$. We  consider the strong solution $u_1(r,t)$ with such initial datum and then perform the transformation
\begin{equation}
u_k(x,t):={\mathcal T}_k u(x,t)= k^{N} u_1(kx, k^{N(p-2)+sp}t)
\end{equation}
for every $k>1$. We want to let $k\to \infty$ in the end. We will apply the continuous rescaling transformation  and study the rescaled flow \eqref{eq.rescflow} (with $a=0$). First, a lemma.

\begin{lemma} If $v_1$ is the rescaled function form $u_1$ and $v_k$ from $u_k$, then
\begin{equation*}
v_k(y,\tau)=  v_1(y, \tau + h)\,,  \quad h=\log(k).
\end{equation*}
This means that the transformation ${\mathcal T}_k$ on the original semigroup becomes a forward time  shift in the rescaled semigroup
\begin{equation}
{\mathcal S}_h v(t)= v(t+h), \quad h=\log(k).
\end{equation}
\end{lemma}

\noindent {\sl Proof.} We have
$$
v_k(y,\tau)=(t+1)^{\alpha}u_k(y(t+1)^{\beta},t)=
k^{N}(t+1)^{\alpha}u(ky(t+1)^{\beta},k^{1/\beta}t),
$$
$$
v_k(y,\tau)=e^{\tau\alpha}u_k(ye^{\tau\beta},e^\tau)=
k^{N}e^{\tau\alpha}u_1(kye^{\beta\tau},k^{1/\beta}e^\tau), )
$$
where $t=e^{\tau}$, $\tau>-\infty$. Put $k=e^{\beta h}$ so that $ke^{\beta\tau}=e^{\beta(\tau +h)}$.
Then 
$$
v_k(y,\tau)=e^{(\tau+h) \alpha}u_1(ye^{\beta(\tau+h)},e^{\tau+h}),
$$
But the inverse transformation  gives $u_1(x,t)=t^{-N\beta}v_1(y',\tau')$, $y'=x\,t^{-\beta},$ $ \tau'=\log(t)\,,$
so that
$$
v_k(y,\tau)= e^{(\tau+h) \alpha} u_1(ye^{\beta(\tau+h)},e^{(\tau+h)})=
e^{(\tau+h-\tau') \alpha} v_1(y'e^{\beta(\tau+h-\tau')}, \tau'(e^h t))
$$
Putting $\tau'=\tau+h$, we get $v_k(y,\tau)=  v_1(y, \tau + h)\,.$ \quad \qed

\noindent $\bullet$  We may pass to the limit in the original family $\{u_k(x,t)\}_k$ or in the rescaled family $\{v_h(y,\tau)\}_h$. The latter is more convenient since it is just the orbit $v_1(\tau)$ and its forward translations.
We will work in finite time intervals $0<t_1\le t\le t_2$, that means $-\infty<  \tau_1\le \tau\le   \tau_2$. From the boundedness estimates we know that both families are bounded and more precisely, the $v$-sequence has a  uniform bound that does not depend on $h$. The family is also uniformly bounded in $L^1(\ren)$. We also have uniform estimates on $v_t$ in $L^\infty_t(L^2_x)$ and $v $ in $L^\infty_t(W^{s,p}_x)$ (Hint: transform the ones for $u_k$). Using the Aubin-Lions compactness results as presented in Simon's \cite{Simon86}, the orbit forms a relatively compact subset of $L^1(\ren)\cap L^2(\ren) $. Therefore, we can pass to the limit $h\to\infty$ and get a limit $V$ with  strong convergence in $\ren\times [t_1,t_2]$.

\noindent The limit $V(y,\tau)$ is a nonnegative solution of the rescaled equation \eqref{eq.resc} for $\tau\ge \tau_i$ with some initial value at $\tau_1$. It satisfies the same bounds as before so it is strong solution in all $L^q$ spaces for $\tau>t_1=-C$. The function is radially decreasing and symmetric in space for all times. The mass is conserved thanks to the uniform tail decay.

\noindent $\bullet$ Going back to the original variables by inverting  transformation \eqref{eq.rescflow}, we get
$$
U(x,t)=t^{-\alpha} V(x\,t^{-\beta}, \log t)\\.
$$
This a strong solution of the original equation \eqref{frplap.eq} that has all the aforementioned properties.
Let check the initial trace. Using the barrier for $u=u_1$ and its decay \eqref{decay} there is a $C>0$ such that
$$
   u_1(x,t)\le C |x|^{-(N+sp)}(t+1)^{-sp\beta}\,.
$$
for all and $t>0$ and $x\ge C(t+1)^{\beta}$. It follows that
$$
   u_k(x,t)\le C k^{N} |kx|^{-(N+sp)}(k^{1/\beta}t+1)^{sp\beta}= C  |x|^{-(N+sp)}(t+k^{-1/\beta})^{sp\beta} \,.
$$
for all $x\ge C(t+k^{-1/\beta})^{\beta}$. In the limit this means that $U(x,t)\le C  |x|^{-(N+sp)}\,t^{sp\beta}$, thus $U$ has a Dirac delta as initial data.
The self-similar solution is constructed. \qed

\medskip

\noindent {\bf Remarks. } 1) It is easy to see that set of self-similar solutions $\{ U_M \} $ is invariant under the mass preserving scaling ${\mathcal T}_k$. In other terms,  the corresponding set of v-solutions $\{ V_M \} $  is invariant under the time translations ${\mathcal S}_h$. This has an important consequence; if we prove uniqueness of the general fundamental solution as constructed in this section, then it  would imply self-similarity because it would imply that such $V$ is stationary in time, hence $U$ is self-similar.
We will not pursue that path in this paper.

 2) Whenever the given total mass is negative, $M<0$, the fundamental solution is obtained by just putting
$U_M(x,t)=-U_{-M}(x,t)$.

\noindent 3) Any  fundamental solution must be radial and decreasing. Use approximation of $\delta$ by $u(\cdot,t)$, with $t$ very small and cut to small support and bounded data and use Proposition \ref{Alek}.
\normalcolor


\subsection{The fundamental self-similar solution}

Since we  did not address the question  of uniqueness in the previous section, we study next the issue of existence of such a self-similar solution. It will be obtained by a method that in a first step  proves existence of periodic $v$ solutions.

\begin{theorem}\label{thm.exfs2}
 There is a fundamental solution of Problem \eqref{frplap.eq}-\eqref{init.data}  with the properties of Theorem \ref{thm.exfs} that is also self-similar. Moreover, the self-similar fundamental solution is unique.
The profile $F$ is a nonnegative and radial $C^1$ function that is nonincreasing along the radius, is positive everywhere and goes to zero at spatial infinity like $O(r^{-(N+sp)})$.
\end{theorem}

\noindent {\sl Proof of existence.} (i) Let $X=L^1(\ren)$. We consider the subset $K\subset X$ consisting of all nonnegative radial functions $\phi$, decreasing along the radial variable,  with mass $\|\phi\|_1\le 1$, and bounded above by one $G$ as in the barrier construction of Theorem \ref{thm.barr}. The  set $K$ is a non-empty, convex, closed and bounded subset with respect to the norm of the Banach space $X$.  Moreover, we have proved that
\begin{equation}\label{bdd.below}
S_\tau\phi(y)\ge \zeta(y)   \quad \mbox{ for } \ y\in \ren, \ \tau>0,
\end{equation}
for a function $\zeta\ge 0$ as in Lemma \ref{lem.pos}.

(ii) Next, we consider the solution of the $v$-equation \eqref{eq.resc} starting at $\tau=0$ with data $v(y,0)=\phi(y)\in K$, and consider the semigroup map $S_h:  X\to X$ defined by $S_h(\phi)=v(\cdot,h)$.   According to our analysis, the set of images $S_h(K)$ satisfies $S_h(K)\subset K$. Moreover, it is relatively compact in $X$. It follows from the Schauder Fixed Point Theorem that there exists at least fixed point $\phi_{h}\in K$, i.\,e., $S_h(\phi_h)=\phi_h$.

Iterating the equality, we  get periodicity for the orbit $v_h(y,\tau)$ starting at $\tau=0$:  for all integers $k\ge 1$ we have
$$
v_h(y,\tau+ kh))=v_h(y,\tau) \quad \forall \tau>0.
$$
 By estimate \eqref{bdd.below} $v_h\ge \zeta$, hence it is not the zero function. Also $v_h\le G$ and it has a certain smoothness.

\medskip

(iii) We now consider the obtained collection of data $\phi_h$ producing periodic such orbits $v_h$ of period $h>0$ and contained in $K$. We may pass to the limit along a subsequence of the dyadic sequence $h_n=2^{-n}$ as $n\to \infty$ and thus find a limit solution $\widehat v$ defined for all $\tau\ge 0$
and starting in $K$, such that the equality
$$
\widehat v(y,\tau+ k2^{-n})=\widehat v(y,\tau) \quad \forall \tau>0
$$
holds for infinitely many $n$'s and all integers $k\ge 1$. By continuity of the orbit in $X$, $\widehat v$ must be stationary in time. Again we conclude that $G(y)\ge \widehat v(y)\ge \zeta$. Going back to the original variables, it means that the corresponding function $\widehat u(x,t)$ is a self-similar solution of equation \eqref{frplap.eq}. Hence, its initial data
must be a non-zero Dirac mass. If it does not have unit mass, at least it has a positive mass. Then, we may use the rescaling \eqref{scal.trn2} to get a self-similar fundamental solution with mass just 1. \qed

The fixed point idea can be found in the literature on asymptotic problems. We mention  Escobedo and Mischler \cite{EMR05} in the study of the equations of coagulation and fragmentation.

\medskip

\noindent {\sl Proof of uniqueness of self-similar profile.} We know that any self-similar profile $F$ is bounded, radially symmetric and non increasing.
 We know that $0\le F\le C$, that $F\le Cr^{-(N+sp)}$.
We prove regularity for the profile by using the regularity of the equation. We recall that
$U_t(x,1)=-\nabla \cdot(xF) $ is bounded, so that $F$ is a $C^1$ function for $r>0$.

 The main step is to use mass difference analysis, since this is a strict Lyapunov functional,
hence we arrive at a contradiction when two self-similar profiles meet. This is an argument
taken from the book \cite{Vazpme07}. It goes as follows: We take two profiles $F_1$ and $F_2$ and assume the same mass $\int F_1\, dx=\int F_2\, dx=1$. If $F_1$ is not $F_2$ they must intersect and then $\int (F_1-F_2)_+dx=C$ is not zero. By self-similarity it must be constant. But we have proved that whenever $C>0$ at one time, it must be a decreasing quantity in time. \qed

\begin{figure}[h!]
\includegraphics[width=\textwidth]{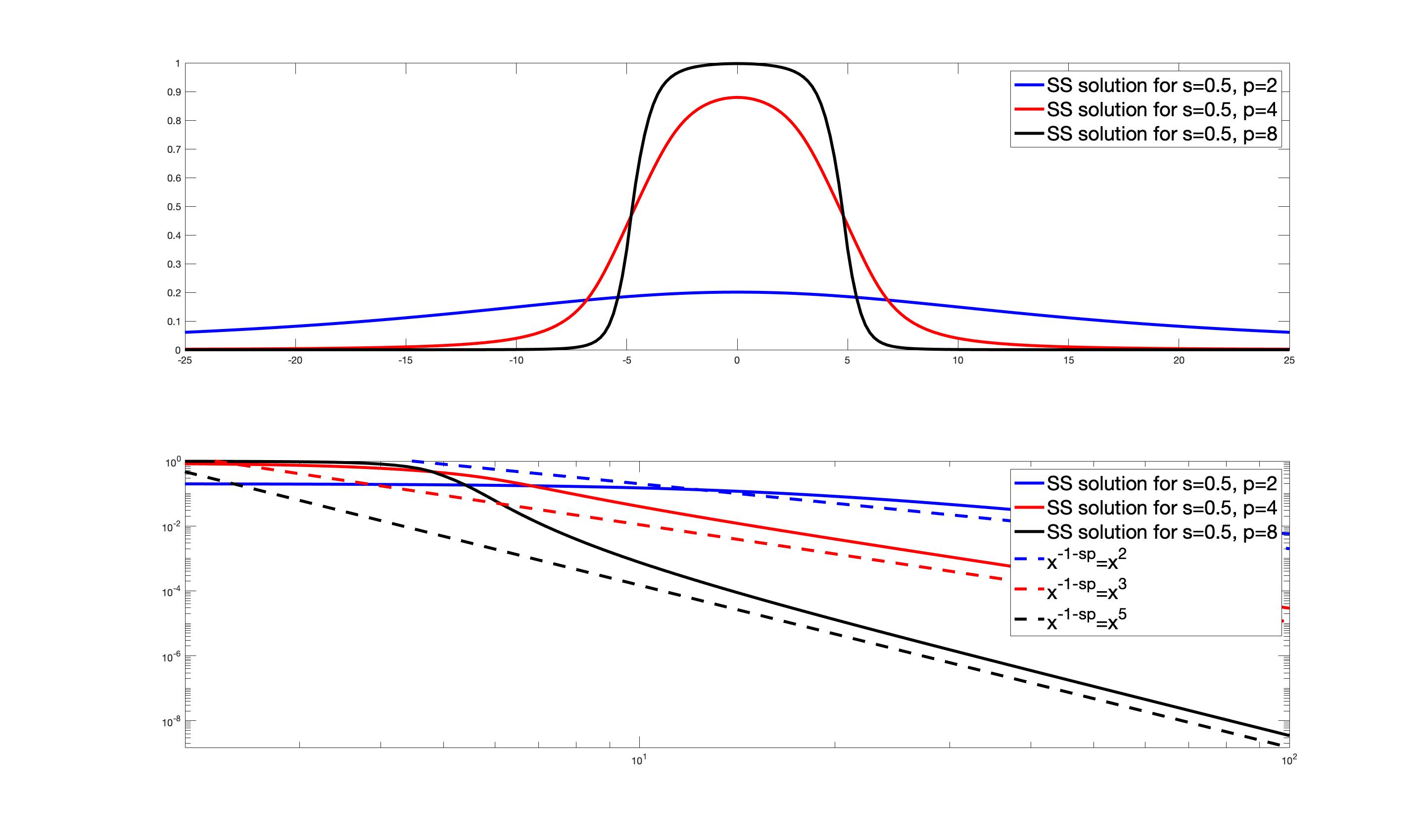}
\caption{Self-similar fundamental solutions for different $p$, with $s=0.5$. }
\label{fig:Self-sim graphics}
\end{figure}

\begin{figure}[h!]
\includegraphics[width=\textwidth]{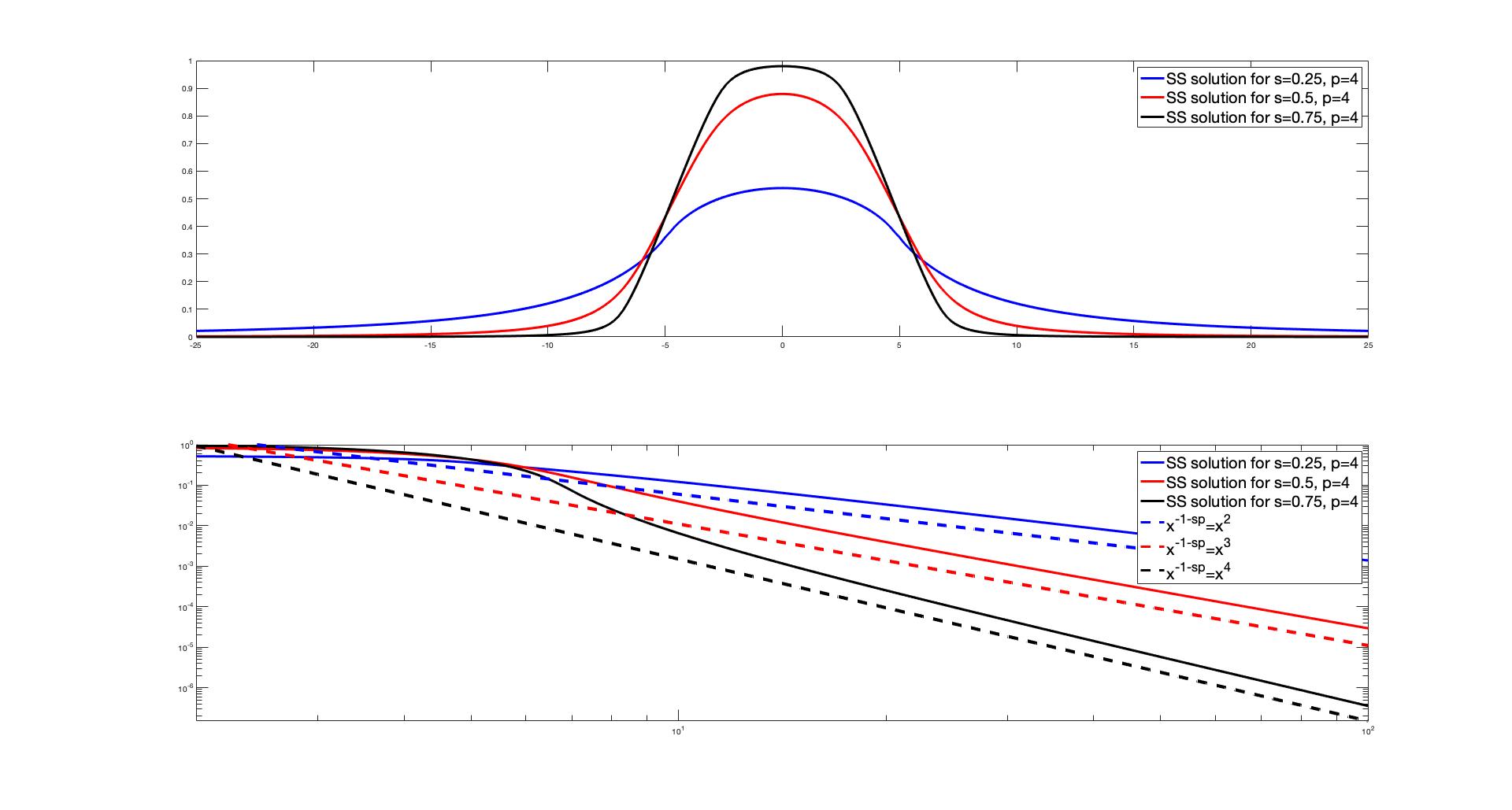}
\caption{Self-similar fundamental solutions for different $s$, with $p=4$.
 }
\label{fig:Self-sim graphics}
\end{figure}

\medskip

\noindent{\bf Computed graphics. }
 Figures 1 and 2 below show the self-similar fundamental solutions for different values of $s$ and $p$. The profiles are computed in dimension $N=1$. The second picture in each  figure shows clearly the predicted decay with exponent $1+sp$ using the logarithmic scale. Also to be remarked the flat behaviour of the profile near the origin for large values of $p$. The numerical treatment is due to F. del Teso.

\section{Positivity and precise tail behaviour}\label{sec.pos.tail}

The fact that solutions of the EFPL equation with nonnegative initial data become immediately positive for all times $t>0$ in the whole space has been proved  in Subsection \ref{ssec.pos}.
Here we will give a more quantitative version of this positivity result.
We recall that in the limit case $s = 1$, with $p > 2$ fixed, we get the standard $p$-Laplacian equation, where positivity at infinity for all nonnegative solutions is false due to the property of
finite propagation. This explains that some special characteristic of fractional
diffusion must play a role to make positivity true.

Our analysis will allow us to obtain the minimum  behavior of nonnegative solutions when $|x|\rightarrow \infty$, more precisely their rate of space decay, for small times $t>0$. This will imply the precise decay rate of the profile of the fundamental solution. Our new idea is to obtain a lower bound that matches the spatial behaviour of the upper barrier, as established in Section \ref{sec.barr}.

\begin{theorem}\label{low.par.est1}
Let $0<s<1$ and $p>2$. Let $\underline u(x,t)$ be a solution of Problem \eqref{frplap.eq} with initial data $u_0(x)\ge 0$ such that $u_0(x)\ge 2$ in the ball $B_2(0)$. Then there is a time
$t_1>0$ and a constant $c>0$ such that
\begin{equation}
u(x,t)\ge c\,t\,|x|^{-(N+sp)}
\end{equation}
if $|x|\ge 2$  and $0<t<t_1$.
\end{theorem}

We will use a comparison argument based on the following construction.

\begin{lemma} There is a smooth, positive and radial function $G_1(r) $ in $\ren$ such that

(i) $G_1(r)\le 1$ everywhere, and $G_1(r)=cr^{-(N+sp)}$ for all $r>2$

(ii) ${\mathcal L}_{s,p}G_1$ is bounded and ${\mathcal L}_{s,p}G_1(r)\approx - r^{-(N+sp)}$ for all $r\ge R>2$.
\end{lemma}

\noindent {\sl Proof of the Lemma.} We define $G_1$ by  specifying it in three different regions.  For $r\le 1$ we put $G=1$.
For $r>2$ we put $G(r)=cr^{-(N+sp)}$ as indicated, with a small constant $0<c<c_0$ that will change in the application, so we must pay attention to it. In the intermediate region we choose a smooth and radially decreasing function that matches the values at $r=1$ and $r=2$ with $C^1$ agreement.

It is then easy from the theory to prove that ${\mathcal L}_{s,p}G_1$ is bounded on any ball, so we only have to worry about the behaviour at infinity, more precisely for $r\gg 2$. In order to analyze that situation  we point out that, according to formula \eqref{frplap.op}, \ ${\mathcal L}_{s,p}G(x)$ is an integral with contributions from the variable $y$ in different regions. We will show that the contribution from the ball $B_1(0)$ is the largest. Indeed, we have for $r>2$
$$
I(r)=-\int_{B_1} \frac{1-G(x)}{|x-y|^{N+sp}}\,dy \le -(1-\ve)|B_1| \ (r/2)^{-(N+sp)},
$$
that does not depend on the small parameter $c$. The other contributions depend on $c$ and can be made small with respect to $I(r)$ for all $r>3$, see details in Section \ref{sec.barr}. \qed

\noindent {\sl Proof of the Theorem.} (i) We modify function $G_1$ to introduce a linear dependence on time in the outer region. We take a smooth cutoff function $\eta$ lying between 0 and 1 such that
$\eta(x)=1$ for $|x|\le 1$ and $\eta(x)=0$ for $|x|\le 2$ and put
\begin{equation}\label{subsol1}
U(x,t)=\eta(x)\,G_1(x)+ (1-\eta(x))ct \,r^{-(N+sp)}.
\end{equation}

\noindent (ii) We want to prove that this function satisfies the subsolution condition
\begin{equation}\label{subsol2}
U_t+{\mathcal L}_{s,p} U <0
\end{equation}
in an outer region $\{r>R\}$ and for an interval of times $0<t<t_*(c)$ if $c$ is small enough. Now for $R>2$ we have
$$
U_t= cr^{-(N+sp)}>0
$$
On the other hand, the proof of the Lemma shows that in that region
$$
{\mathcal L}_{s,p} U \le -Cr^{-(N+sp)},
$$
as long as we can disregard the contributions from outside $B_1$, and this is true if $tc$ is mall enough. The conclusion \eqref{subsol2} follows.

\noindent (iii) We now to the comparison step between $u$ and $U$ in a space-time domain of the form
$Q=\{(x,t): \ |x|\ge 3, \ 0<t<t_1\}$.  By comparison we may consider some smaller initial data $u_0$, such that $0\le u_0(x)\le 2$ and $u_0(x)=  2$ in the ball of radius 3. Moreover, $u_0$ is smooth. By previous results of this paper we know that $u(x,t) \in C^\alpha (\ren \times [0,T])$ and $u(x,t)>0$ for all $x\in\ren$ and $t>0$. We have that  $u(x,t)\ge 1$ in a ball of radius $2<R<3$ for all small times $0<t<t_0$.

We already have the necessary inequalities for the equation inside that domain. We must check the initial and lateral outside conditions.

As for initial conditions we know that $U(x,0)=0$ for all $|x|\ge 2$, while $u_0\ge 0$ everywhere.

Regarding comparison for $|x|\le 3$ we know that $U(x,t)\le 1$ at all points (while $t$ is small)
while the continuity of the solution $u$ and its initial data  imply that $u(x,t)\ge 2-\ve$ for all $|x|\ge 2$ and $0<t<t_*$.

\noindent (iv)  Now we only need to integrate by parts the difference of the two equations with multiplier $(U-u)_+$ to get the conclusion that $(u-U)_+$ must be zero a.e. in $Q$. Note that both functions belong to $L^2(\ren)\le L^\infty(\ren)$ uniformly in $t$. Since $U\le u$ in the set $\Omega=\{|x|\ge 3\}$ for $0<t<t_0$, we get for all those times
$$
\frac{d}{dt}\int_{\Omega} (U-u)_+^2\,dx= 2\int_{\Omega} (U-u)_+(U_t-u_t)\,dx=
2\int_{\ren} (U-u)_+(U_t-u_t)\,dx=I.
$$
But that integral is easily estimated
$$
I=-2\int_{\ren} ({\mathcal L}_{s,p}U-{\mathcal L}_{s,p}u)\,(U-u)_+\,dx\le 0
$$
by $T$-accretivity (better do the direct computation, see above the computation of the evolution of the $L^2$  norm of the difference of two solutions). Since  $(U-u)_+=0$ in $\Omega$ for $t=0$, we get the desired conclusion:
$$
u(x,t)\ge U(x,t)\ge ct\,r^{-(N+sp)}
 $$
 if  $r\ge 3$ and  $t<t_0$. \qed


\subsection{Application to the self-similar solution} We consider the fundamental solution after a time displacement:
$$
u_1(x,t)= (t+1)^{-\alpha} F_M(|x|\,(t+1)^{-\beta}),
$$
that satisfies the assumptions of Theorem \ref{low.par.est1} if $M>0$ is large enough. We conclude that
$$
F_M(r)\ge C_1\,r^{-(N+sp)} \quad \mbox{for all large } \ r.
$$
By scaling, the same is true for $M=1$ with a different constant. Together with the upper bound from Theorem \ref{thm.barr}, the last assertion of Theorem \ref{thm.ssfs} is proved.

For possible future reference, let us state the tail behaviour of the fundamental solution $U_M(x,t)=t^{-\alpha}F_M(|x|\,t^{-\beta})$. Let us choose $M>0$.

 \begin{corollary} \label{cor.ob.fs} On every outer region of the form $\{(x,t): \ |x|\ge C t^{\beta}, \ C>0\}$ we have constants $0<C_1<C_2$ such that
\begin{equation}\label{sigma}
C_1 \,M^{\sigma}|x|^{-(N+sp)}t^{sp\beta}\le U_M(x,t)\le C_2\, M^{\sigma}|x|^{-(N+sp)}t^{sp\beta},
\end{equation}
where $\sigma=1+(p-2)sp\beta$.
 \end{corollary}

\medskip

\noindent {\bf Remark.} Positivity estimates related to the ones in this section have been obtained for the fractional porous medium equation in \cite{VazBar2014, StanVa14, VolVa14}. Other forms of positivity estimates  were developed in \cite{BV06} for the Fast Diffusion Equation, and in the fractional case in \cite{BV14}.

\section{Asymptotic Behaviour}\label{sec.ab}

We establish here the asymptotic behaviour of finite mass solutions, reflected in Theorem \ref{thm.ab1}.
We may assume that $M>0$ and the case $M<0$ can be reduced to positive mass by changing the sign of the solution.
We comment on $M=0$ below.

(i) We prove first the $L^1$ convergence.  By scaling we may also assume that $M=1$.  The proof   relies on the previous results plus the existence of a strict Lyapunov functional,  that happens to be
\begin{equation}
J(u_1,u_2;t):=\int (u_1(x,t)-u_2(x,t))_+\,dx
\end{equation}
where $u_1$ and $u_2$ are two solutions with finite mass.

\begin{lemma}
Let $u_1$ and $u_2$ are two solutions with finite mass. Then, $J(u_1,u_2;t)$ is strictly decreasing in time
unless the solutions are ordered.
\end{lemma}

\noindent {\sl Proof. } By previous analysis, Section \ref{sec.diff},
we know that
\begin{equation}
\frac{d}{dt}J(u_1,u_2;t)=-\iint_D \left| |u_1(x:y,t)|^{p-2}u_1(x:y,t) -
|u_2(x:y,t)|^{p-2}u_2(x:y,t)\right|\,d\mu(x,y)\,,
\end{equation}
with notation as in \eqref{diff.est.1b}. In particular, the set $D\subset \re^{2N}$ contains the points where
$$
(u_1(x,t)-u_2(x,t))\,(u_1(y,t)-u_2(y,t))<0.
$$
Now, in order to $dJ/dt$ to vanish at a time $t_0>0$ we need $u_1(x:y,t)=u_2(x:y,t)|$ on $D$, i.\,e., $u_1(x,t)- u_2(x,t)=u_1(y,t)- u_2(y,t)$. But this is incompatible with the definition $D$, so $D$ must be empty, hence $u_1$ and $u_2$ must be ordered at time $t$. This implies that they have the same property for $t>t_0$. \qed

\medskip

\noindent {\sl  Proof of Theorem \ref{thm.ab1} continued.} It is convenient to  consider the $v$ version of both solutions, namely $v_1$ and $V_M$. We can show that $v(y,\tau+ n_k)$ converges strongly in $L^1(\ren)$, along a subsequence $n_k\to \infty$, towards a new solution $w_1$ of the $v$-equation. Under our assumptions $w_1$ is a fundamental solution. On the other hand, $V_M$ is stationary.

We know from the Lemma that $J(v_1, V_M;t)$ is strictly decreasing in time, unless $v_1(t)=V_M$ for all large $t$, in which case we are done. If this is not the case, we continue as follows. By monotonicity there is a limit
$$
\lim_{t\to\infty} J(u_1, U_M;t)=\lim_{\tau\to\infty} J(v_1, V_M;\tau)=C\ge 0\\.
$$
We want to prove that $C=0$, which implies our result. If the limit is not zero, we consider the evolution of the new solution $w_1$ together with $V_M$. We have
$$
 J(w_1, V_M;t_0)=\lim_{\tau \to \infty}J(v_1, V_M; t_0+\tau)=C\,,
$$
i.e., is constant for all $t_0>0$, which means that $w_1=V_M$  by equality of mass \normalcolor and the lemma. By uniqueness of the limit, we get convergence along the whole half line $t>0$ instead of a sequence of times.

For general data $u_0\in L^1(\ren)$, $M>0$, we use approximation.

Finally, in the case $M=0$ we just bound our solution from above and below by solutions of mess $\ve$ and $-\ve$ resp,,
apply the Theorem and pass to the limit $\ve\to 0$.

\noindent (ii) {\sl  Proof of  convergence in uniform norm, formula  \eqref{lim.ab.inf}.} We return to the proof of the previous step and discover that the bounded sequence $v(y,\tau+ n_k)$ is locally relatively compact in the set of continuous functions in $\ren\times (\tau_1, \tau_2) $ thanks to the results on H{\"o}lder continuity of \cite{BLS2019} as commented in Subsection \ref{ssec.bddness}, once they are translated to the $v$-equation. Hence, it converges locally to the same limit as before, but now in uniform norm. In order to get global convergence we need to control the tails at infinity. We use the following argument: a sequence of space functions $v(\cdot, \tau)$ that is uniformly bounded near infinity in $L^1$ (thanks to the convergence to $V_M$) and is also uniformly H{\"o}lder continuous must also be also uniformly small in $L^\infty$. This implies that the previous uniform convergence was not only local but global in space. Using the correspondence \eqref{eq.rescflow},  we get the convergence of the $u(t)$ with factor $t^{\alpha}$. This part of the theorem is proved.   \qed

\section{Two-sided global bounds. Global Harnack }\label{sec.gharn}
The uniform convergence of the previous section implies that
$u(x,t)/U_M(x,t)\to 1$ as $t\to1\infty$ uniformly on sets of the form $\{|x|\le ct^{\beta}\}$. But it does not say anything about the relative error on the far away region, i.e. for the so-called tail behaviour.

We can contribute to that issue using the positivity analysis of Section \ref{sec.pos.tail}. We obtain a two-sided global estimate, assuming that the initial data are bounded, nonnegative and compactly supported. The result applies to all positive times and says that the relative quotient $u(x,t)/U_M(x,t)$ stays bounded for $t\ge \tau>0$.

\begin{theorem}\label{thm.GH} Let $u$ the semigroup solution corresponding to initial data $u_0\in L^\infty{\ren}$, $u_0\ge 0$, $u_0\ne 0$, supported in a ball of radius $R$. For every $\tau>0$ there exist constants $M_1, M_2>0$ and delay $ c_2>0$ such that
\begin{equation}
U_{M_1}(x,t)\le u(x,t)\le U_{M_2}(x,t+c_2) \qquad \mbox{ for all \ } \ x\in \ren, \ t\ge \tau.
\end{equation}
The constants  $M_1, M_2$, and $c_2$ may depend on $\tau$. Moreover, if $M(u_0)=\int_{\ren} u_0\,dx$, then $M_1\le M(u_0)\le M_2$.
\end{theorem}

\medskip

\noindent {\sl Proof}. (i) Let us begin by the upper bound that is an easy consequence of the barrier estimate of Section \ref{sec.barr}, in particular Theorem \ref{thm.barr}. Indeed, the function $G$ mentioned there  is comparable with the self-similar profile $F_1$, hence smaller than $F_{M_2}$ for some $M_2>1$. This estimate is valid even for $\tau=0$, with $c_2=1$ and $M_2$ large enough. It holds then for every $t>0$ by comparison.

\noindent (ii) For the lower bound we need to take $\tau>0$ and use the results of this section. By comparison, translations and rescaling we my assume that $u_0$ is as in Theorem \ref{low.par.est1}. We also assume that is radially decreasing. Therefore, given a time  $\tau>0$ small enough we have the estimate \ $ u(x,\tau)\ge c\,\tau\,|x|^{-(N+sp)} $ \
for all $|x|\ge R > 2$. On the other hand,
$$
U_{M_1}(x,c_1)= c_1^{-\alpha}F_{M_1}(|x|\,c_1^{-\beta})\le C \min\{ c_1^{-\alpha}\, M_1^{sp\beta}, \,M_1^{\sigma} c_1^{sp\beta}|x|^{-(N+sp)}\},
$$
for every $|x|\ge 0$. We have used formula \eqref{sigma}. We conclude that for given $c_1>0$ there exists $M_1$ small enough such that
$$
u(x,\tau)\ge U_{M_1}(x,c_1).
$$
We may put now $\tau=c_1$. By comparison the lower bound is true for all $t\ge \tau$.

(iii)  In view of the previous asymptotic results we have $M_1<M(u_0)<M_2$. Just check the values at $x=0$ for large $t$. \qed

This kind of two-sided bound by the fundamental solution is usually called a Global Harnack Inequality and is frequent in nonlinear diffusion problems with fast diffusion. See  applications to the fast diffusion equation in \cite{Vascppme, CrVa03, BV06}, and a very recent one in \cite{Simonov}. It is not true for equations with slow diffusion and free boundaries. There are a number of references for fractional parabolic equations like  \cite{BSV15,  BSV17}, even in the so-called slow range, like the present paper.  We will not mention the large literature on elliptic problems or problems in bounded domains.  

Let us write in a clear way our conclusion about the size of the spatial tails.

\begin{corollary} For every solution $u(x,t) $ and as in the previous theorem, and for every $t\ge 1$ there are positive constants $C_1,C_2>0$ such that
\begin{equation}
C_1< u(x,t)\,|x|^{N+sp}\,t^{-sp\beta}<C_2
\end{equation}
on the outer set  $|x|\ge t^{\beta}$.
\end{corollary}

\section{Source-type solution in a bounded domain}\label{sec.st-sol}

We can derive from the previous study the existence of source-type  solutions for the problem posed in a bounded domain with zero Dirichlet outside conditions.  They take a Dirac delta as initial data but we do not call them fundamental because they do not play such a key role in the theory.

\begin{theorem}\label{thm.st-sol} There exists a solution of the Dirichlet problem for equation \eqref{frplap.eq}
posed in a bounded domain $\Omega\subset \ren$ with initial data a Dirac delta located at an interior  point, $x_0\in \Omega$, and zero Dirichlet data outside $\Omega$. For $t\ge \tau>0$, it is a bounded strong solution of the equation as described in \cite{Vaz16}.
\end{theorem}

\noindent {\sl Proof.}  (i) For convenience, we assume in the first step that $\Omega$ is the ball radius 1 centered at 0 and $x_0=0$. We may also assume that $M=0$. Existence and uniqueness of solutions for the Cauchy-Dirichlet has been established in \cite{Vaz16} and other references, and an ordered semigroup of contractions is generated in all $L^q$ spaces, $1\le q<\infty$. Further estimates and  regularity are obtained, but beware of the long-time behaviour that is completely different. Here a question of small time behaviour is of concern, and luckily there is great similarity in that issue.

(ii) The existence of solutions of the approximate problems with data $u_{0n}\ge 0$ that converge to a Dirac delta does not offer any difficulty. Passing to the limit we easily obtain a solution $U(x,t)$ of the Cauchy-Dirichlet problem in $B_1$, using the a priori estimates and known compactness. The only important missing point is justifying that the initial data are taken. We recall that mass is not conserved in time for the Cauchy-Dirichlet problem in a bounded domain.

In order to solve the pending issue, it will be enough to show that the mass of the limit solution $U(\cdot,t)$ tends to 1 as $t\to 0$.
We want to prove that for an approximating sequence of functions $u_{0n}\ge 0$, $\int_{B_1} u_{0n}(x)\,dx=1$ and $u_{0n}(x)\to \delta(x)$ weakly,
then for every $\epsilon>0$ there is an $n_0$ and a $t_0$ such that
\begin{equation}\label{eps.claim}
\int_{B_1}u_{n}(x,t)\,dx> 1-\epsilon, \forall n\ge n_0, \ 0<t<t_0.
\end{equation}

We take the same initial data $u_{0n}$ as an approximating sequence for the problem in $\ren$ and in this way we show that the corresponding solutions that we now call $u_n^\infty(x,t)$ converge to the self-similar fundamental solution that we call $U^\infty(x,t)$, and we have described in previous sections. By comparison we have
$$
u_n(x,t)\le u_n^\infty(x,t), \quad U(x,t)\le U^\infty(x,t).
$$

(iii) The novelty comes next. The following lemma provides a proof of the needed estimate \eqref{eps.claim}. We will also assume that the initial data $u_{0n}$ are a sequence of rescalings of an initial $u_{01}$ that is nonnegative, smooth, bounded and supported in a small ball $B_\delta(0)$.

\begin{lemma} Under the previous assumptions, for every $\ve>0$ there are $n_0$ and $\tau$ such that for $n\ge n_0$  the following inequality holds
$$
u_n^\infty(x,t)-\ve\le u_n(x,t) \ \mbox{ in \ } \ B_1(0)\times (0,\tau).
$$
Therefore, $ U(x,t)\ge U^\infty(x,t)-\ve$ in $B_1(0)\times (0,\tau)$.
\end{lemma}

\noindent {\sl Proof.} We first claim that \ $\widetilde u(x,t)=u_n^\infty(x,t)-\ve$ is a solution of the same equation \eqref{frplap.eq}  posed in the context of the space $X_\ve$ obtained from $L^1(\ren)$ by subjecting  all functions to a downward shift. This is due to the fact the operator in invariant under vertical shifts.  After the shift, the initial data are lower that before in $B_1$. In the exterior of the ball, $|x|\ge 1$,  $u_n(x,t)$ is extended but zero, while we can check that for large $n$
\begin{equation}\label{ineq.ext}
u_n^\infty(x,t)-\ve\le 0 \quad \mbox{for all} \ |x|\ge 1, \ \mbox{for  \ } \ 0<t<\tau,
\end{equation}
thanks to the a priori estimates on the decay of the solutions. Admitting this fact for the moment, we may now use comparison of the solutions in the ball to conclude that
$u_n^\infty(x,t)-\ve\le u_n(x,t) $ { in } $ B_1(0)\times (0,\tau)$ as desired, and this implies \eqref{eps.claim}.

In order to prove \eqref{ineq.ext} we use the a priori estimate for all the sequence $u_n$ in terms of the barrier as stated at the end of Section \ref{sec.barr}
$$
u_n(x,t)\le C |x|^{-(N+sp)}(t+a)^{sp\beta}\,.
$$
This constant depends on the initial data. We need $u_{on}$ to be below the barrier at $t=0$ and for that need that for $n$ large and putting $|x|= \delta/n$ we have
$$
c_1n^N\le C(\delta/n)^{-(N+sp)}a^{sp\beta},
$$
i.\,e., $C\ge c_1 \delta^{-(N+sp)} n^{-sp}a^{-sp\beta}$ near infinity.  \color{blue} We conclude that we can fix a uniform $C$ at for $n\le n_0$.  \normalcolor We go back to the outer comparison. We need
$$
C |x|^{-(N+sp)}(t+a)^{sp\beta}\le \ve
 $$
for $|x|\ge 1$ and $0<t<\tau$. This holds if $ C (\tau+a)^{sp\beta}\le \ve$. \qed


\subsection{Other domains}

(i) We consider first the case of balls $B_R$ of radius $R>0$. Given some initial data $u_0\in  L^1(B_1)$ we can solve the Cauchy-Dirichlet problem in $B_1$ to obtain a function $u(x,t)=S_t(u_0)$, where $S_t$ is the semigroup generated by the equation in $B_1$. Likewise, we denote the semigroup in $B_R $ by $S_t^{R}$, and the semigroup in $\ren$ by $\overline S_t$.

It is easy to see that the scaling \ $T_R u(x,t) = R^{-N }u(x(R,t/R^{1/\beta})$ \
generates a function $u^{R}=T_Ru$ that solves the same Cauchy-Dirichlet problem in $B_R$. Moreover,
$$
u^{R}(x,0):=T_R u_0(x) = R^{-N }u(x/R)
$$
is a rescaling of  $u_0$ that is defined for all $x\in R$. Mass in conserved (at corresponding times). We have $S_t^R(u^R(0))= T_R S_t(u_0)$.
The transformation can be inverted using $(T_R)^{-1}=T_{1/R}$. It is clear that $T_R$ transforms a source-type solution in $B_1$ into a  source-type solution in $B_R$. Besides, the Maximum Principle implies that for all $u_0\in L^1(\ren)$, $u_o\ge 0$ we have
$$
S_t(u_0)\le S_t^{R}(u_0)\le \overline S_t(u_0),
$$
A similar order applies to  fundamental solutions.

(ii) For other domains $\Omega\subset \ren$ we use comparison with balls to make sure that the usual approximate solutions so not lose the initial trace when passing to the limit.
More precisely, after translation we may assume that $0\in \Omega$ and that $B_{R_1}(0)\subset \Omega\subset B_{R_2}(0)$.
In this way the existence of a source-type solution in $\Omega$ is proved. We leave the details to the reader.

\section{Limit cases}\label{sec.limitcases}

In the paper we have considered  all fractional exponents in the range $0<s<1$ and nonlinear exponents $p>2$. The limit cases are interesting as examples of continuity with the dependence on parameters. We will make here a brief sketch of important facts.

\medskip

\noindent $\bullet$ {\sl Limit $p\to 2$.} The limit of the $(s,p)$-semigroup as $p\to 2$ for fixed $s$ offers only minor difficulties.   Also the passage to the limit in the self-similar solutions gives the well-known profiles of the fractional linear heat equation. These profiles decay like $O(|x|^{N+2s} )$ as $|x|\to\infty$, cf. \cite{Blumenthal-Getoor}, see also \cite{BSV17} and its references. The linear self-similar solutions are also asymptotic attractors, as proved in \cite{Va18}, where convergence rates are obtained. The  limit  $p\to 2$ can also be checked computationally with minor difficulty.

\noindent $\bullet$  {\sl Limit $s\to 1$.} It is well known by experts that the operator ${\mathcal L}_{s,p}$ must be conveniently renormalized by a constant including the factor $1-s$,  cf. \cite{BBM02, IshiiN, MazRT}, in order to converge to the standard $p$-Laplacian as $s\to 1$. It is then rather easy to prove that, as $s\to 1$ for fixed $p> 2$, we obtain the semigroup corresponding to the standard $p$-Laplace operator, already mentioned in the introduction.

In particular, we can pass to the limit in the self-similar solutions of Theorem \ref{thm.ssfs} and check that the self-similar profile $F_{s,p}(r)$ converges to the  profile $F_{1,p}(r)$, given by the well known Barenblatt profile
$$
F(r)=\left(C- kr^{\frac{p}{p-1}}\right)_+^{\frac{p-1}{p-2}},
$$
cf. \cite{VazSmooth}, formula 11.8. The decay exponents $\alpha(s,p)$ and $\beta(s,p)$ also converge. Notice that for $0<s<1$ the self-similar profiles $F_{s,p}(r)$ are positive with tails at infinity, but the limit Barenblatt profile, $F_{1,p}(r)$, has compact support.

Full details should be provided elsewhere. A similar study of continuous dependence with respect to parameters has been done in full detail in the case of the Fractional Porous Medium Equation in \cite{DPQRV2}.
\nc

\medskip

\noindent $\bullet$  {\sl Limit $s\to 0$.}  This case offered a very nice surprise to the author. We want to take the limit $s\to 0$ in the fundamental solutions $F_s(y)=F(y;s,p,M)$ that  have been constructed and described above. We look at the equation satisfied by $F_s$:
$$
{\mathcal L}_{s,p} F= \beta \,\nabla\cdot (yF)\,,
$$
and pass to the limit $s\to 0$. With  a proper scaling ${\mathcal L}_{s,p}$ tends to the identity, $\alpha\to 1/(p-2)$ and $\beta\to 1/(N(p-2))$. We get in a formal way
an equation for any limit profile $F(r)$,  which turns out to be a simple ODE:
\begin{equation}\label{ODE}
 NF+ rF'= N(p-2)F^{p-1}\,.
\end{equation}
As limit of the self-similar profiles $F_s$, the profile $F$  for $s=0$ must be nonnegative with $F'(r)\le 0$. An analysis of the ODE shows that $F(r)$ must live in the rectangle
$$
0<r<\infty, \quad 0\le F(r)\le F_* =(p-2)^{-1/(p-2)}
$$
where $F_*>0$ is the value that corresponds to a constant solution. All other positive nonincreasing  solutions of the ODE \eqref{ODE} start at $F(0)=F_*$ and decrease to the value $F(+\infty)=0$
with an asymptotic estimate of the form
$$
F(r;C)\sim C\,r^{-N},
$$
hence they have infinite mass. And we are lucky since \eqref{ODE} is a Riccatti equation with exact solutions
$$
F(r)=\frac1{((p-2) + C\,r^{N(p-2)})^{1/(p-2)}}\,,
$$
which gives in the $(x,t)$ variables
$$
U(x,t)=((p-2)t + C\,|x|^{N(p-2)})^{-1/(p-2)}\,.
$$
We point out that this is a solution of the formal limit of the parabolic equation which is
$$
u_t=-u^{p-1}.
$$
For the simplest situation, we may take $p=3$ and then
$$
F(r)=(1+Cr^{N})^{-1}, \quad U(x,t)=(t+C|x|^{N})^{-1}.
$$
Note that this limit solutions are not integrable, in other words, they have infinite mass.
This subsection is an announcement of new facts. Full details of the limit process should be provided elsewhere.

\section{Final comments}

We begin with a technical appendix. Then, we complement the information on related work given in the introduction with some historical comments of the main topic of the paper. The last subsection may interest the curious reader.

\subsection{\bf Calculation of the $s$-$p$-Laplacian for $C^2$ functions}\label{ssec.ta}

This a technical reminder for the reader. We want to prove that, when applied to a function
$u\in C^2(\ren)$ with  bounded norms, the $s$-$p$-Laplacian has a well-defined value for every $x\in\ren$, and moreover, it  is a continuous function. We assume that $0<s<1$ and $p>2$. By definition
$$
{\mathcal L}_{s,p}(u)(x)=\int \frac{|u(x)-u(x-y)|^{p-2}(u(x)-u(x-y))}{|y|^{N+sp}}\,dy=
$$
$$
\frac12\int (|u(x)-u(x-y)|^{p-2}(u(x)-u(x-y))+ |u(x)-u(x-y)|^{p-2}(u(x)-u(x-y)))\,\frac{dy}
{|y|^{N+sp}}.
$$
Now we use the inequality for $p>2$
$$
||a|^{p-2}a\pm|b|^{p-2} b|\le C(p)||a|^{p-2}+|b|^{p-2}||a\pm b|
$$
Apply this formula with $a=u(x)-u(x-y)$ and $b=u(x)-u(x+y)$ to get an estimate for the integrand:
\begin{gather*}
|(u(x)-u(x-y)|^{p-2}(u(x)-u(x-y))+ |u(x)-u(x+y)|^{p-2}(u(x)-u(x+y))| \le\\
C(p)(|u(x)-u(x-y)|^{p-2}+|u(x)-u(x+y)|^{p-2})\,|2u(x)-u(x+y)-u(x-y)|\le\\
2C(p)|Du(x')y|^{p-2} |D^2u(x'')y^2|\le 2C(p)\|Du\|_\infty^{p-2}\|D^2 u\|_\infty |y|^{p}.
\end{gather*}
This proves that the integral is uniformly convergent for functions $u\in C^2(\ren)$ with bounded $L^\infty$ norms.  The resulting integral is a continuous function of $x$. Moreover, we get the interpolation formula
\begin{equation}
|{\mathcal L}_{s,p}(u)(x)| \le C_1\|u\|_{\infty}^{p-1}+ C_2\|Du\|_{\infty}^{p-2}\|D^2u\|_{\infty}\,.
\end{equation}
Hint: split the integral into the domains $|y|\le 1$ and $|y|\ge 1$.

This type of calculation is also used in \cite{DTGV20} in the study of different representations of the $s$-$p$-Laplacian. For a more delicate calculation valid for all $p\in (1,\infty)$, see \cite{KKL}, Section 3.
 \qed


\subsection{Fundamental solutions in nonlinear diffusion} \label{sec.cfs}
The importance of the Gaussian fundamental solution in the classical heat equation is well-known in the mathematics literature and needs no reminder, \cite{EvansPDE, Widder}. In the linear fractional case $p=2$ with $0<s<1$,  the fundamental solution of the fractional heat equation is also known thanks to  Blumental and Getoor \cite{Blumenthal-Getoor} who studied it in 1960.   In such a  case the fundamental solution also allows to construct the class of all nonnegative solutions of the Cauchy problem in the whole space by using the representation formula, see the theory of \cite{BSV17} where an optimal class of data is considered and well-posedness shown.

In the case of  nonlinear problems, the importance of fundamental solutions  has been proved in numerous examples, even if, contrary to what happens for linear equations, representation formulas for general solutions in terms of such a special solution are not available. Their interest lies mainly in the description of the asymptotic behaviour as $t\to\infty$  of general solutions. The fundamental solution is well-known in the standard $p$-Laplacian case, $p> 2$, $s=1$. Its existence comes from \cite{Bar}, hence the name Barenblatt solution, and its uniqueness was established in  \cite{KV88}, see also \cite{Laur98}. For the standard porous medium equation the situation is well-known, see the historical comment in the monograph \cite{Vazpme07}. A  recent example for nonlinear fractional equations is given by the fundamental solution of the fractional porous medium equation constructed by the author in \cite{VazBar2014}. For the so-called  porous medium equation with fractional potential pressure the fundamental solution was first constructed in   \cite{BKM} and \cite{CaVa11}, and the asymptotic behaviour was established in the latter reference. In all cases the application to the asymptotic behaviour as $t\to\infty$ is carried through, and convergence of a general class of finite-mass solutions to the corresponding fundamental solution is proved.

On the other hand, for the problem posed in  a bounded domain the special solution that is relevant  concerning the asymptotic behaviour as $t\to\infty$ is the separate-variables solution called the {\sl friendly giant}. It was constructed for our equation by the author in  \cite{Vaz16}.

\subsection{Other comments and extensions}

\noindent $\bullet$ We have proved uniqueness of the self-similar fundamental solution. The uniqueness of the general fundamental solution is a delicate issue that we did not settle here.

\noindent $\bullet$ The  exact tail behaviour  of the fundamental solution  may be improved.
The numerical computations suggest a finer decay expression $F(x)\sim  C\, |x|^{N+sp}$.\nc

\noindent $\bullet$ Existence of solutions for measures as initial data should be investigated.  This is related to the question of initial traces.

\noindent $\bullet$ The question of rates of convergence for the result \eqref{lim.ab.L1} of Theorem  \ref{thm.ab1} has not been considered. This issue has been addressed for many other models of nonlinear diffusion. It is solved for many of them, but well known cases remain open.

\noindent $\bullet$ We did not consider the case where $1<p<2$, which has its own features. \nc For $p$ close to 2 there exists a fundamental solution that explains the asymptotic behaviour, much as done here. This property is well known for the standard $p$-Laplacian equation with an explicit formula, cf. \cite{VazSmooth}, formula (11.8). Likewise,  there is a critical exponent for our equation when the self-similarity exponents blow up, i.e., for $p_c=2N/(N+s)$. For $p<p_c$  such a fundamental solution does not exist. There is extinction in finite time, as proved in \cite{BS2020}.

\noindent $\bullet$  In the existence theory we can consider wider classes of initial data, possibly growing at infinity. Optimal classes are known in the linear fractional equation (case $p=2$), \cite{BSV17}, and in the standard $p$-Laplacian equation (case $s=1$), cf. \cite{DiBeHerr89}. Of course, the asymptotic behaviour will not be the same.

\noindent $\bullet$ Another interesting issue is the presence of a right-hand side in the equation, maybe in the form of lower-order teems. There are some works, see e.g. 
\cite{Teng2019} and its references.

\noindent $\bullet$ We have considered a nonlinear equation of fractional type with nonlinearity $\Phi(u)=|u|^{p-2}u$, and we have used the fact that $\Phi$ is a power, hence homogeneous, in a number of tools. We wonder how much of the theory holds for more general monotone nonlinearities $\Phi$.


\vskip .6cm

\noindent {\textbf{\large \sc Acknowledgments.}} Author partially funded by Projects MTM2014-52240-P and   PGC2018-098440-B-I00 (Spain). Partially performed as an Honorary Professor at Univ. Complutense de Madrid. We thank M. Bonforte for information of ongoing  work with A. Salort \cite{BS2020}, that covers a number of topics of this paper. Thus, it treats the case $p<2$, but does not treat the self-similar solutions or asymptotic behaviour. We thank J. L. Rodrigo for pointing out reference \cite{EMR05} and E. Lindgren for information on his work \cite{BLS2019}. The numerical study leading to the figures is due to F. del Teso.
\normalcolor

\medskip

{\small
\bibliographystyle{amsplain}

}

\

\noindent {\sc Address:}

\noindent Juan Luis V\'azquez. Departamento de Matem\'{a}ticas, Universidad
Aut\'{o}noma de Madrid,\\ Campus de Cantoblanco, 28049 Madrid, Spain.  \\
e-mail address:~\texttt{juanluis.vazquez@uam.es}

\medskip

\noindent {\bf 2020 Mathematics Subject Classification.}
  	35K55,  	
   	35K65,   	
    35R11,   	
    35A08,   	
    35B40.   	

\noindent {\bf Keywords: } Nonlinear parabolic equations, $p$-Laplacian operator, fractional operators, fundamental solutions,  asymptotic behaviour.

%

\end{document}